\documentclass[3p]{elsarticle}
\usepackage[utf8]{inputenc}
\usepackage{lineno}
\usepackage{amsfonts}
\usepackage{graphicx}
\usepackage{epstopdf}
\usepackage{algorithmic}
\usepackage{mathtools}
\usepackage{amsmath}
\usepackage{enumerate}
\usepackage{upgreek}
\usepackage{subcaption, caption}
\usepackage[T1]{sansmath} 
\usepackage{lmodern,bm} 
\usepackage{pdflscape}
\usepackage[pdfencoding=auto]{hyperref}
\usepackage{booktabs}
\usepackage{tikz} 
\usetikzlibrary{shapes.multipart}
\usetikzlibrary{matrix}
\usetikzlibrary{positioning}
\usepackage{ragged2e}
\usepackage{verbatim} 
\usepackage{soul}

\usepackage{calc}
\newsavebox\CBox
\newcommand\hcancel[2][0.5pt]{%
  \ifmmode\sbox\CBox{$#2$}\else\sbox\CBox{#2}\fi%
  \makebox[0pt][l]{\usebox\CBox}%
  \rule[0.5\ht\CBox-#1/2]{\wd\CBox}{#1}}

\DeclarePairedDelimiter\floor{\lfloor}{\rfloor}

\DeclareMathAlphabet\bfcal{OMS}{cmsy}{b}{n}
\newcommand{\chinf}{\chi}

\journal{Journal of Computational Physics}
\begin{document}

\begin{frontmatter}

\title{Accelerated Calder\'on preconditioning for Maxwell transmission problems}

\author[ucl]{Antigoni Kleanthous}
\ead{antigoni.kleanthous.12@ucl.ac.uk}
\author[ucl]{Timo Betcke}
\author[ucl]{David P. Hewett}
\author[AI]{Paul Escapil-Inchausp\'e}
\author[AI]{Carlos Jerez-Hanckes}
\author[met,hertfordshire]{Anthony J. Baran}

\address[ucl]{University College London, Department of Mathematics, London, UK}
\address[AI]{Universidad Adolfo Ib\'a\~nez, Faculty of Engineering and Sciences, Santiago, Chile}
\address[met]{Met Office, Exeter, UK}
\address[hertfordshire]{University of Hertfordshire, School of Physics, Astronomy, and Mathematics, Hertfordshire, UK}

\begin{abstract}
We investigate a range of techniques for the acceleration of Calder\'on (operator) preconditioning in the context of boundary integral equation methods for electromagnetic transmission problems. 
Our objective is to mitigate as far as possible the high computational cost of the barycentrically-refined meshes necessary for the stable discretisation of operator products. 
Our focus is on the well-known PMCHWT formulation, but the techniques we introduce can be applied generically. 
By using barycentric meshes only for the preconditioner and not for the original boundary integral operator, we achieve significant reductions in computational cost by (i) using ``reduced'' Calder\'on preconditioners obtained by discarding constituent boundary integral operators that are not essential for regularisation, and (ii) adopting a ``bi-parametric'' approach \cite{escapil2019fast, ESCAPILINCHAUSPE2021220} in which we use a lower quality (cheaper) $\mathcal{H}$-matrix assembly routine for the preconditioner than for the original operator, including a novel approach of discarding far-field interactions in the preconditioner. 
Using the boundary element software Bempp (www.bempp.com), we compare the performance of different combinations of these techniques in the context of scattering by multiple dielectric particles. 
Applying our accelerated implementation to 3D electromagnetic scattering by an aggregate consisting of 8 monomer ice crystals of overall diameter 1cm at 664GHz leads to a 99\% reduction in memory cost and at least a 75\% reduction in total computation time compared to a non-accelerated implementation.
\end{abstract}

\end{frontmatter}

\section{Introduction}
\label{sec:introduction}

Boundary Integral Equations (BIEs) are a powerful tool in the simulation of electromagnetic scattering. 
In this approach, a Maxwell scattering problem is reduced to a BIE of the form
\begin{align}
    \label{eqn:BIE}
    \bm{\mathcal{A}} \mathbf{u}^s = \mathbf{f},
\end{align}
where the solution $\mathbf{u}^s$ involves certain traces of the unknown scattered electric and/or magnetic fields on the boundary of the scatterer, the data $\mathbf{f}$ depends on appropriate boundary traces of the known incident field, and $\bm{\mathcal{A}}$ is a bounded linear operator (in appropriate function spaces), involving one or both of the electric and magnetic boundary integral operators, which we denote in this paper by $\mathcal{S}$ and $\mathcal{C}$, respectively. 

Typically, upon standard Galerkin discretisation \eqref{eqn:BIE} produces a highly ill-conditioned linear system requiring a large number of iterations when an iterative solver, such as GMRES, 
is used. At the operator level this can be attributed to the hypersingular nature of $\mathcal{S}$ and the compact nature of $\mathcal{C}$ (at least on smooth domains, and again, in appropriate function spaces, discussed below), which at the discrete level can lead to eigenvalues accumulating at both infinity and zero.

Calder\'on preconditioning is an operator preconditioning approach \cite{hiptmair2006operator}, based on the fact that $\mathcal{S}$ and $\mathcal{C}$ are, in a certain sense, self-regularising, satisfying 
\begin{align}
    \label{eqn:calderon_identities}\mathcal{S}^2 = -\frac{1}{4}\mathcal{I} + \mathcal{C}^2, \qquad \mathcal{C} \mathcal{S} + \mathcal{S} \mathcal{C} = 0.
\end{align}
These identities suggest the possibility of defining a preconditioning operator 
$\bm{\mathcal{P}}$, involving a suitable combination of the operators $\mathcal{S}$ and $\mathcal{C}$, so that the discretised form of the preconditioned system
\begin{align}
    \label{eqn:preconditioned_BIE}\bm{\mathcal{PA}} \mathbf{u} =  \bm{\mathcal{P}}\mathbf{f},
\end{align}
can be solved more efficiently than that of \eqref{eqn:BIE}. 
The overall computational cost of solving \eqref{eqn:preconditioned_BIE} is strongly affected by both 
\begin{enumerate}[(i)]
\item \label{i} the choice of preconditioning operator $\bm{\mathcal{P}}$; and 
\item \label{ii} the choice of discretisation for the operator product $\bm{\mathcal{P}}\bm{\mathcal{A}}$. 
\end{enumerate}
Regarding \eqref{i}, one obviously wants $\bm{\mathcal{P}}$ to be an effective regularizer of $\bm{\mathcal{A}}$, in the sense of reducing the number of GMRES iterations. But one would also like $\bm{\mathcal{P}}$ to involve as few instances of the operators $\mathcal{S}$ and $\mathcal{C}$ as possible, so as to minimise the cost of each iteration. 
Regarding \eqref{ii}, one would like to use a discretisation for the product $\bm{\mathcal{P}}\bm{\mathcal{A}}$ that is as cheap as possible, while maintaining sufficient solution accuracy. 
However, to stably discretise $\bm{\mathcal{P}}\bm{\mathcal{A}}$ using boundary elements one needs to work with both a primal mesh and a dual mesh, the latter defined in terms of a barycentric refinement with six times as many elements as the primal mesh. For high frequency problems, when one is already using a fine mesh to capture the oscillatory solution, this requirement of barycentric refinement can be prohibitive. For this reason, the acceleration of Calder\'on preconditioning techniques for electromagnetic scattering is currently an active research area in numerical analysis and scientific computing.

Our focus in the current paper is on Maxwell transmission problems modelling scattering by homogeneous dielectric bodies, one application of which is in the simulation of scattering of electromagnetic radiation by ice crystals in cirrus clouds, important in climate and weather modelling (see \cite{baran2009review,baran2012single,liou2016light,kleanthous2019calderon} and the references therein). 
We confine our attention to the Poggio-Miller-Chang-Harrington-Wu-Tsai (PMCHWT) formulation of \cite{poggio1970integral, wu1977scattering, mautz1977electromagnetic,harrington1989boundary}, but emphasise that our ideas could also be applied to other BIE formulations of transmission problems such as M\"uller \cite{muller2013foundations} formulations. 

In the context of the PMCHWT formulation, regarding point \eqref{i} above, the classical choice of Calder\'on preconditioner is $\bfcal{P}=\bfcal{A}$ \cite{yan2010comparative, cools2011calderon, niino2012calderon, yla2015challenges} (for other related work see \cite{contopanagos2002well, christiansen2002preconditioner, antoine2008integral, bagci2009calderon, boubendir2015integral}).
It was shown in \cite{kleanthous2019calderon} that for multi-particle configurations it is more efficient to use $\bm{\mathcal{P}} = \bm{\mathcal{D}}$, where $\bm{\mathcal{D}}$ is the block-diagonal component of $\bm{\mathcal{A}}$ whose diagonal blocks are the PMCHWT operators for the individual scatterers. Indeed, for this choice approximately the same number of GMRES iterations are required compared to the classical choice $\bm{\mathcal{P}} = \bm{\mathcal{A}}$, but the computational cost per iteration is lower. 
We also refer to \cite{thierry2014remark} where the idea of the reduced Calder\'on preconditioner $\bm{\mathcal{P}} = \bm{\mathcal{D}}$ is applied to multi-particle acoustic scattering.
Our contribution to \eqref{i} in the current paper is to show that further cost savings can be obtained at high frequencies by further reducing the complexity of $\bfcal{P}$, so as to include only a quarter of the operators appearing in each block of $\bfcal{D}$.

Our contribution to \eqref{ii} is to demonstrate that for the PMCHWT formulation the cost of discretising the operator product $\bm{\mathcal{P}}\bm{\mathcal{A}}$ can be reduced, without sacrificing solution accuracy, by using a cheaper (poorer quality) matrix assembly routine for $\bm{\mathcal{P}}$ than for $\bm{\mathcal{A}}$. 
Our approach is based on an adaptation of the ``bi-parametric'' $\mathcal{H}$-matrix approach proposed in \cite{escapil2019fast} in the context of the electric field integral equation (EFIE) for perfectly conducting  scatterers. For the Helmholtz case, we refer to \cite{FJH20} and for a general theoretical framework see \cite{ESCAPILINCHAUSPE2021220}. Specifically, following \cite{escapil2019fast} we adopt $\mathcal{H}$-matrix approximations for the assembly and storage of the discrete versions of $\bm{\mathcal{P}}$ and $\bm{\mathcal{A}}$ in which we
\begin{itemize}
    \item use a weaker target tolerance for $\bm{\mathcal{P}}$ than for $\bm{\mathcal{A}}$ in the approximation of admissible blocks by low-rank approximations; and 
    \item use lower order quadrature rules for $\bm{\mathcal{P}}$ than for $\bm{\mathcal{A}}$.
\end{itemize}
In addition, going beyond the framework of \cite{escapil2019fast}, for $\bfcal{P}$ (but not $\bfcal{A}$) we assemble admissible blocks only if they correspond to sufficiently `near-field' interactions, discarding far-field contributions in the preconditioner.

To get the most benefit from these acceleration techniques we found it advantageous to choose a combination of discrete trial and test spaces such that $\bfcal{P}$ is discretised using the more expensive dual mesh and $\bfcal{A}$ is discretised using the cheaper primal mesh. 
The high cost of the barycentric refinement is then mitigated by our use of a reduced-complexity preconditioner with a low-cost assembly routine. 

Our numerical results show that the combination of a reduced Calder\'on preconditioner with a bi-parametric implementation results in significant reductions in computational cost, especially at high frequencies. For the highest frequency problem considered here, relating to electromagnetic scattering by a realistic complex aggregate of dielectric particles used to represent an ice crystal aggregate of the type which occurs naturally in cirrus clouds, the total computation time was reduced by 75\% while the memory cost was reduced by 99\%, compared to the non-accelerated approach.

The structure of the paper is as follows. In \S\ref{sct:scattering_problem} we describe the scattering problem and the PMCHWT BIE formulation. In \S\ref{sct:calderon_construction} we discuss Calder\'on preconditioners and their discretisations, and outline our proposed acceleration techniques. In \S\ref{sct:benchmarks} we report a selection of results from a detailed numerical study into the performance of our techniques on a benchmark problem, and in \S\ref{sec:Applications} we consider applications to scattering by a realistic complex ice aggregate particle.
Finally, in \S\ref{sct:conclusion} we present some conclusions. 
All of our numerical experiments were performed using the boundary element software library Bempp (www.bempp.com) \cite{smigaj2015solving}, Version 3.3.5. \footnote{Example notebooks are available at \url{https://github.com/ankleanthous/Accelerated_Calderon}.}

\section{The scattering problem and the PMCHWT BIE formulation}\label{sct:scattering_problem}
We consider 3D time-harmonic electromagnetic scattering (with $\mathrm{e}^{-\mathrm{i} \omega t}$ time-dependence) by $M$ disjoint isotropic homogeneous (but not necessarily identical) dielectric scatterers occupying bounded domains $\Omega^i_m\subset \mathbb{R}^3$, $m=1,\ldots,M$, with boundaries $\Gamma_m = \partial \Omega^i_m$, in a homogeneous exterior medium $\Omega^e = \mathbb{R}^3 \backslash \cup_{m=1}^M\overline{\Omega^i_m}$.
For $m=1,\ldots, M$, we denote the (constant) electric permittivity and magnetic permeability of the scatterer $\Omega^i_m$ by $\epsilon_m$ and $\mu_m$, respectively, and the resulting wavenumber by $k_m = \omega \sqrt{\mu_m \epsilon_m}$. We denote the corresponding parameters of the exterior medium $\Omega^e$ by $\epsilon_e$, $\mu_e$, and $k_e = \omega \sqrt{\mu_e \epsilon_e}$. For $m=1,\ldots, M$ we define the (in general complex) refractive index of the scatterer $\Omega^i_m$ to be the ratio $n_m=k_m/k_e$. 

In the scattering problem, an incident field $(\mathbf{E}^{inc}, \mathbf{H}^{inc})$ in $\Omega^e$ creates interior fields $(\mathbf{E}^i_m, \mathbf{H}^i_m)$ in $\Omega^i_m$, for $m = 1, \ldots, M$, and a scattered field $(\mathbf{E}^s, \mathbf{H}^s)$ in $\Omega^e$, which is assumed to satisfy the Silver-M\"uller radiation condition, and which combines with the incident field to give the total exterior field
\begin{align}
\label{eqn:TotalFields}
    \mathbf{E}^e & = \mathbf{E}^s + \mathbf{E}^{inc}, \qquad
    \mathbf{H}^e  = \mathbf{H}^s + \mathbf{H}^{inc}, \qquad\text{in }\Omega^e.
\end{align}
The interior and exterior fields satisfy the Maxwell equations 
\begin{alignat}{4}
\label{eqn:maxwell_i}\nabla \times \mathbf{E}^i_m &= \mathrm{i} \omega \mu_m \mathbf{H}^i_m, 
&\qquad
\nabla \times \mathbf{H}^i_m &= - \mathrm{i} \omega \epsilon_m \mathbf{E}^i_m, &&\qquad\text{in }\Omega^i_m,\,m=1,\ldots,M,\\
\label{eqn:maxwell_e}\nabla \times \mathbf{E}^e &= \mathrm{i} \omega \mu_e \mathbf{H}^e, 
&\qquad
\nabla \times \mathbf{H}^e &= - \mathrm{i} \omega \epsilon_e \mathbf{E}^e,&&\qquad\text{in }\Omega^e,
\end{alignat}
together with the transmission boundary conditions
\begin{alignat}{3}
\label{bc}\mathbf{E}^i_m (\mathbf{x}) \times \mathbf{n}_m &= \mathbf{E}^e (\mathbf{x}) \times \mathbf{n}_m, 
&\qquad
\mathbf{H}^i_m (\mathbf{x}) \times \mathbf{n}_m &= \mathbf{H}^e (\mathbf{x}) \times \mathbf{n}_m, &&\qquad  \mathbf{x} \in \Gamma_m, \ m=1,\ldots,M,
\end{alignat}
where $\mathbf{n}_m$ is the unit normal vector on $\Gamma_m$ pointing into $\Omega^e$. Since the magnetic fields can subsequently be recovered from \eqref{eqn:maxwell_i}-\eqref{eqn:maxwell_e}, it is sufficient to solve for the electric fields alone, which satisfy
\begin{align}
\label{eqn:curlcurl_i}
\nabla \times \nabla \times \mathbf{E}^i_m - k_m^2 \mathbf{E}^i_m = \mathbf{0}, 
&\qquad\text{in }\Omega^i_m,\,m=1,\ldots,M,\\
\label{eqn:curlcurl_e}
\nabla \times \nabla \times \mathbf{E}^e -k_e^2 \mathbf{E}^e= \mathbf{0}, 
&\qquad\text{in }\Omega^e.
\end{align}

Before reviewing the PMCHWT BIE formulation we fix some notation. 
For a bounded Lipschitz open set $\Omega$ with boundary $\Gamma = \partial \Omega$ and outward unit normal vector $\mathbf{n}$, and for some wavenumber $k$, we define the Dirichlet and Neumann traces of vector fields $\mathbf{P}^\pm$ defined in the exterior ($+$) and interior ($-$) of $\Omega$ by
\begin{alignat}{4}
\gamma_D^{\pm} \mathbf{P}^\pm (\mathbf{x}) &:=& \mathbf{P}^\pm(\mathbf{x}) \times \mathbf{n},
&\qquad 
\gamma_N^{\pm} \mathbf{P}^\pm(\mathbf{x}) &:=& \frac{1}{\mathrm{i}k} \gamma_D^{\pm} \left( \nabla \times  \mathbf{P}^\pm(\mathbf{x}) \right), &\qquad \mathbf{x} \in \Gamma,
\end{alignat}
and the electric and magnetic potential operators acting on a boundary vector field $\mathbf{v}$ by
\begin{align}
\mathcal{E} \mathbf{v}(\mathbf{x}) &:= \mathrm{i}k \int _\Gamma \mathbf{v}(\mathbf{y}) G (\mathbf{x}, \mathbf{y}) \mathrm{d} s (\mathbf{y}) - \frac{1}{\mathrm{i} k}  \nabla_\mathbf{x} \int _\Gamma \nabla_\mathbf{y} \cdot \mathbf{v}(\mathbf{y}) G (\mathbf{x}, \mathbf{y}) \mathrm{d} s (\mathbf{y}),\\
\mathcal{H} \mathbf{v}(\mathbf{x}) &:= \nabla_\mathbf{x} \times \int _\Gamma \mathbf{v}(\mathbf{y}) G (\mathbf{x}, \mathbf{y}) \mathrm{d} s (\mathbf{y}),
\end{align}
where $G (\mathbf{x}, \mathbf{y}) = \frac{\exp (\mathrm{i} k | \mathbf{x} - \mathbf{y}|)}{4 \pi | \mathbf{x} - \mathbf{y}|}$ is the fundamental solution of the 3D Helmholtz equation. 
The associated electric and magnetic boundary integral operators, which also act on boundary vector fields, can be defined by $\mathcal{S}:= \frac{1}{2}(\gamma_D^+ + \gamma_D^-)\mathcal{E}$ and $\mathcal{C}:= \frac{1}{2}(\gamma_D^+ + \gamma_D^-)\mathcal{H}$, and satisfy the relations
\begin{alignat}{2}
\label{BIOrelns}
\mathcal{S} &= \gamma_{D}^\pm \mathcal{E} = - \gamma_{N}^\pm \mathcal{H}, &\qquad
\mathcal{C} &= \gamma_{N}^\pm \mathcal{E} \pm \frac{1}{2} \mathcal{I} = \gamma_{D}^\pm \mathcal{H} \pm \frac{1}{2} \mathcal{I},
\end{alignat}
where $\mathcal{I}$ is the identity operator. These formal statements can be made rigorous with $\mathcal{S}$ and $\mathcal{C}$ defining bounded linear endomorphisms on the space $\mathbf{H}_\times^{-\frac{1}{2}} (\textnormal{div}_\Gamma, \Gamma)$ of tangential vector fields (currents) on $\Gamma$ that, along with their divergences, have Sobolev regularity $-\frac{1}{2}$ (for details see \cite{buffa2003galerkin}). 

Having introduced this notation, we now return to the scattering problem. By the Stratton-Chu formulae \cite{kirsch2015mathematical} the electric fields satisfying \eqref{eqn:TotalFields}-\eqref{eqn:maxwell_e} can be represented 
as
\begin{alignat}{3}
\label{eqn:StrattonChu_int} \mathcal{H}^i_m (\gamma_{D,m}^- \mathbf{E}^i_m) &+ \mathcal{E}^i_m (\gamma_{N,m}^- \mathbf{E}^i_m) = 
\begin{cases}
\mathbf{E}^i_m(\mathbf{x}), & \mathbf{x} \in \Omega^i_m, \\
\mathbf{0}, & \mathbf{x} \not \in \overline{\Omega^i_m},
\end{cases}\\
\label{eqn:StrattonChu_ext}-\sum_m^M \mathcal{H}^e_m (\gamma_{D,m}^+ \mathbf{E}^s) &- \sum_m^M  \mathcal{E}^e_m (\gamma_{N,m}^+ \mathbf{E}^s) = 
\begin{cases}
\mathbf{E}^s(\mathbf{x}), & \mathbf{x} \in \Omega_e, \\
\mathbf{0}, & \mathbf{x} \not \in \overline{\Omega_e},
\end{cases}
\end{alignat}
where $(\mathcal{E}^i_m,\mathcal{H}^i_m,\gamma_{D,m}^-,\gamma_{N,m}^-)$ are $(\mathcal{E},\mathcal{H},\gamma_{D}^-,\gamma_{N}^-)$ for $\Gamma=\Gamma_m$ and $k=k_m$, and $(\mathcal{E}^e_m,\mathcal{H}^e_m,\gamma_{D,m}^+,\gamma_{N,m}^+)$ are $(\mathcal{E},\mathcal{H},\gamma_{D}^+,\gamma_{N}^+)$ for $\Gamma=\Gamma_m$ and $k=k_e$, for $m=1,\ldots, M$. 
The boundary conditions \eqref{bc} can be re-written 
as 
\begin{align}
    \label{eqn:boundary_conditions_traces} \mathbf{u}^i_m = \mathbf{u}^s_m + \mathbf{u}^{inc}_m,
\end{align}
where
\begin{align}
\mathbf{u}^i_m = \begin{bmatrix}
\gamma_{D,m}^+ \mathbf{E}^i_m \\[6pt]
\dfrac{k_m}{\mu_m} \gamma_{N,m}^+ \mathbf{E}^i_m
\end{bmatrix}, \quad
\mathbf{u}^s_m = \begin{bmatrix}
\gamma_{D,m}^+ \mathbf{E}^s \\[6pt]
\dfrac{k_e}{\mu_e} \gamma_{N,m}^+ \mathbf{E}^s
\end{bmatrix}, \quad
\mathbf{u}^{inc}_m = \begin{bmatrix}
\gamma_{D,m}^+ \mathbf{E}^{inc} \\[6pt]
\dfrac{k_e}{\mu_e} \gamma_{N,m}^+ \mathbf{E}^{inc}
\end{bmatrix}.
\end{align}
By taking Dirichlet and Neumann traces of \eqref{eqn:StrattonChu_int}-\eqref{eqn:StrattonChu_ext} and applying \eqref{eqn:boundary_conditions_traces}, one obtains a set of BIEs satisfied by the unknown boundary traces of $\mathbf{E}^i_m$ and $\mathbf{E}^s$ (for details see \cite{kleanthous2019calderon}), which one can combine in different ways to obtain different BIE formulations. In this paper, we focus on the PMCHWT formulation, which for our multi-particle scattering problem can be written as 
\begin{alignat}{3}
\label{eqn:PMCHWT_multiple_scatterers}\bm{\mathcal{A}}\mathbf{ u}^s = \left(\frac{1}{2}\bm{\mathcal{ {I}}} - \bm{\mathcal{{D}}}^i \right) \mathbf{u}^{inc},
\end{alignat}
where 
\begin{align}
    \label{multiple_operator}\bfcal{A} &= \begin{bmatrix}
    \bfcal{A}_1^e + \bfcal{A}_1^i & \bfcal{A}_{12} & \cdots & \bfcal{A}_{1M} \\[3pt]
    \bfcal{A}_{21} & \ddots & \ddots & \vdots \\[3pt]
    \vdots & \ddots & \ddots & \bfcal{A}_{(M-1)M} \\[3pt]
    \bfcal{A}_{M1} & \cdots & \bfcal{A}_{M(M-1)} & \bfcal{A}_M^e + \bfcal{A}_M^i
    \end{bmatrix}, \quad 
    \bfcal{D}^i = \begin{bmatrix}
    \bfcal{A}_1^i & 0 & \cdots & 0 \\[3pt]
    0 & \ddots & \ddots & \vdots \\[3pt]
    \vdots & \ddots & \ddots & 0 \\[3pt]
    0 & \cdots & 0 & \bfcal{A}_M^i
    \end{bmatrix}, \\
    \label{multiple_operator2} \bfcal{I} &= \begin{bmatrix}
    \bfcal{I}_1 & 0 & \cdots & 0 \\[3pt]
    0 & \ddots & \ddots & \vdots \\[3pt]
    \vdots & \ddots & \ddots & 0 \\[3pt]
    0 & \cdots & 0 & \bfcal{I}_M
    \end{bmatrix}, \quad 
    \mathbf{u}^s = \begin{bmatrix}
    \mathbf{u}^s_1 \\[3pt]
    \mathbf{u}^s_2 \\[3pt]
    \vdots \\[3pt]
    \mathbf{u}^s_M
    \end{bmatrix}, \quad 
    \mathbf{u}^{inc} = \begin{bmatrix}
    \mathbf{u}^{inc}_1 \\[3pt]
    \mathbf{u}^{inc}_2 \\[3pt]
    \vdots \\[3pt]
    \mathbf{u}^{inc}_M
    \end{bmatrix},
\end{align}
with
\begin{align}
\label{AimAem}
\bm{\mathcal{{A}}}^i_m &= \begin{bmatrix}
\mathcal{C}^i_m & \dfrac{\mu_m}{k_m} \mathcal{S}^i_m \\[6pt]
-\dfrac{k_m}{\mu_m} \mathcal{S}^i_m & \mathcal{C}^i_m
\end{bmatrix}, \
\bm{\mathcal{{A}}}^e_m = \begin{bmatrix}
\mathcal{C}^e_m & \dfrac{\mu_e}{k_e} \mathcal{S}^e_m \\[6pt]
-\dfrac{k_e}{\mu_e} \mathcal{S}^e_m & \mathcal{C}^e_m
\end{bmatrix}, \
\bm{\mathcal{{A}}}_{m\ell} = \begin{bmatrix}
    \mathcal{C}^e_{m\ell} & \dfrac{\mu_e}{k_e} \mathcal{S}^e_{m\ell} \\[6pt]
    -\dfrac{k_e}{\mu_e} \mathcal{S}^e_{m\ell} & \mathcal{C}^e_{m\ell}
\end{bmatrix}, \
\bm{\mathcal{I}}_m  = \begin{bmatrix}
\mathcal{I}_m & 0 \\[6pt]
0 & \mathcal{I}_m
\end{bmatrix}.
\end{align}
Here $(\mathcal{C}^i_m,\mathcal{S}^i_m)$ and $(\mathcal{C}^e_m,\mathcal{S}^e_m)$ are $(\mathcal{C},\mathcal{S})$ for $\Gamma=\Gamma_m$ and $k=k_m$ or $k=k_e$, respectively, $(\mathcal{C}^e_{m\ell},\mathcal{S}^e_{m\ell})$ are the analogous ``off-diagonal'' operators, involving integration on $\Gamma_\ell$ and evaluation on $\Gamma_m$, with $k=k_e$, and $\mathcal{I}_m$ is the identity operator on $\mathbf{H}_\times^{-\frac{1}{2}}(\textnormal{div}_{\Gamma_m}, \Gamma_m)$. For details see \cite{kleanthous2019calderon}. 

To obtain a numerical approximation to the solution of the original scattering problem the standard workflow is to first discretise \eqref{eqn:PMCHWT_multiple_scatterers}, using for example a Galerkin discretisation with a piecewise polynomial approximation space, and solve the resulting finite dimensional system of linear equations to approximate the unknown traces $\mathbf{u}^s$. One can then recover the remaining traces using the boundary conditions \eqref{bc}, and finally evaluate the representations \eqref{eqn:StrattonChu_int}-\eqref{eqn:StrattonChu_ext} to obtain the interior and scattered fields. 
However, as alluded to in \S\ref{sec:introduction}, standard Galerkin discretisation of \eqref{eqn:PMCHWT_multiple_scatterers} typically produces ill-conditioned linear systems on which iterative methods such as GMRES converge slowly. 
In the next section, we turn our attention to the design of suitable operator-based preconditioning strategies to remedy this. 

\section{Calder\'on preconditioners and acceleration techniques}\label{sct:calderon_construction}
The PMCHWT operator $\bm{\mathcal{A}}$ defined in \eqref{multiple_operator} is continuous as a mapping $\bm{\mathcal{A}}: \mathbf{X}\rightarrow \mathbf{X}$, where $\mathbf{X} = \oplus_{m=1}^M \mathbf{H}^{-\frac{1}{2}}_\times (\textnormal{div}_{\Gamma_m}, \Gamma_m)^2$. 
The idea of Calder\'on preconditioning is to use the relations \eqref{eqn:calderon_identities} to identify a suitable preconditioning operator $\bm{\mathcal{P}}:\mathbf{X}\to \mathbf{X}$, involving some subset of the operators $\mathcal{C}^i_m$, $\mathcal{C}^e_m$, $\mathcal{C}^e_{m\ell}$, $\mathcal{S}^i_m$, $\mathcal{S}^i_m$ and $\mathcal{S}^e_{m\ell}$, and a suitable discretisation strategy for the operator product $\bm{\mathcal{P}}\bm{\mathcal{A}}$, such that the linear system arising from the discretisation of the preconditioned operator equation 
\begin{alignat}{3}
\label{eqn:PMCHWT_multiple_scatterers_Precond}\bm{\mathcal{P}}\bm{\mathcal{A}}\mathbf{ u}^s = \bm{\mathcal{P}}\left(\frac{1}{2}\bm{\mathcal{ {I}}} - \bm{\mathcal{{D}}}^i \right) \mathbf{u}^{inc}
\end{alignat}
can be solved more efficiently using an iterative method than that arising from \eqref{eqn:PMCHWT_multiple_scatterers}. Our aim in this paper is to study some novel choices of $\bfcal{P}$ and accompanying discretisation strategies and compare their performance on practical problems. 

\subsection{Choice of preconditioning operator $\bfcal{P}$}
\label{sec:Pchoice}
We first consider the choice of preconditioning operator $\bfcal{P}$. 
As mentioned in \S\ref{sec:introduction}, the classical choice for \eqref{eqn:PMCHWT_multiple_scatterers} is $\bfcal{P}=\bfcal{A}$. 
In \cite{kleanthous2019calderon} it was shown that, for multiple scatterers, one can achieve a similar improvement in conditioning with a reduced computational cost by taking $\bfcal{P}=\bfcal{D}$, where 
\begin{align}
    \bfcal{D} = \begin{bmatrix}
    \bfcal{A}_1^e + \bfcal{A}_1^i & 0 & \cdots & 0 \\[3pt]
    0 & \ddots & \ddots & \vdots \\[3pt]
    \vdots & \ddots & \ddots & 0 \\[3pt]
    0 & \cdots & 0 & \bfcal{A}_M^e + \bfcal{A}_M^i
    \end{bmatrix}
\end{align}
is the block diagonal component of $\bfcal{A}$. 
The efficacy of this approach is explained by the fact that the ill-conditioning of $\bfcal{A}$ comes from the diagonal blocks of \eqref{eqn:PMCHWT_multiple_scatterers}, since
\begin{itemize}
\item  the operators $\mathcal{C}^i_m$ and $\mathcal{C}^e_m$ are compact on smooth domains \cite{nedelec2001acoustic}, with eigenvalues clustering around zero;
\item the operators $\mathcal{S}^i_m$ and $\mathcal{S}^e_m$ are sums of compact operators, with eigenvalues clustering around zero, and hypersingular operators, with eigenvalues accumulating at infinity \cite{cools2011calderon}; 
\item and the off-diagonal blocks of $\bm{\mathcal{A}}$ are purely compact, since the kernels of the operators $\mathcal{C}^e_{m\ell}$ and $\mathcal{S}^e_{m\ell}$ are smooth.
\end{itemize}

In this paper, we consider some ``reduced'' versions of the block-diagonal preconditioner, namely 
$\bfcal{P}=\bm{\mathcal{D}}^e$, $\bfcal{P}=\bm{\mathcal{D}}^i$, $\bfcal{P}=\bm{\mathcal{S}}^e$ and $\bfcal{P}=\bm{\mathcal{S}}^i$, where $\bfcal{D}^i$ is as defined in \eqref{multiple_operator} and   
\begin{align}
    \bfcal{D}^e = \begin{bmatrix}
    \bfcal{A}_1^e & 0 & \cdots & 0 \\[3pt]
    0 & \ddots & \ddots & \vdots \\[3pt]
    \vdots & \ddots & \ddots & 0 \\[3pt]
    0 & \cdots & 0 & \bfcal{A}_M^e
    \end{bmatrix}, \quad
    \bfcal{S}^i = \begin{bmatrix}
    \begin{bmatrix}
    0 & \dfrac{\mu_1}{k_1}\mathcal{S}^i_1 \\
    -\dfrac{k_1}{\mu_1}\mathcal{S}^i_1 & 0
    \end{bmatrix} & 0 & \cdots & 0 \\[3pt]
    0 & \ddots & \ddots & \vdots \\[3pt]
    \vdots & \ddots & \ddots & 0 \\[3pt]
    0 & \cdots & 0 & \begin{bmatrix}
    0 & \dfrac{\mu_M}{k_M}\mathcal{S}^i_M \\
    -\dfrac{k_M}{\mu_M}\mathcal{S}^i_1 & 0
    \end{bmatrix}
    \end{bmatrix} 
\end{align} 
\begin{align} 
    \bfcal{S}^e = \begin{bmatrix}
    \begin{bmatrix}
    0 & \dfrac{\mu_e}{k_e}\mathcal{S}^e_1 \\
    -\dfrac{k_e}{\mu_e}\mathcal{S}^e_1 & 0
    \end{bmatrix} & 0 & \cdots & 0 \\[3pt]
    0 & \ddots & \ddots & \vdots \\[3pt]
    \vdots & \ddots & \ddots & 0 \\[3pt]
    0 & \cdots & 0 & \begin{bmatrix}
    0 & \dfrac{\mu_e}{k_e}\mathcal{S}^e_M \\
    -\dfrac{k_e}{\mu_e}\mathcal{S}^e_1 & 0
    \end{bmatrix}
    \end{bmatrix}
\end{align}
Our rationale for studying these reduced versions is that (i) we expect them to still provide some preconditioning effect, because of the presence of the operators $\mathcal{S}^i_m$ and/or $\mathcal{S}^e_m$ in the appropriate positions (note in particular that $\mathcal{S}^i_m$ and $\mathcal{S}^e_m$ can regularise each other since one is just a compact perturbation of the other), and (ii) they involve even fewer boundary integral operators than $\bfcal{D}$, so should require less computational effort to assemble, store and apply within the GMRES iteration. Certainly, for the particular discretisation we use (see the next section for details) the memory cost and assembly times for $\bfcal{S}^i$ or $\bfcal{S}^e$ should be half those of $\bfcal{D}^i$ and $\bfcal{D}^e$, which in turn should be half those of $\bfcal{D}$. If we measure computation time for the linear solve in terms of ``matvecs'', defined as in \cite{kleanthous2019calderon} to mean the time required for a single application of one discretised boundary integral operator $\mathcal{C}^i_m$, $\mathcal{C}^e_m$, $\mathcal{S}^i_m$, $\mathcal{S}^e_m$ etc., then the overall costs for the above formulations (including the initial pre-multiplication of the right hand side where applicable) are given in Table \ref{tab:matvecs}. The choices $\bfcal{P}=\bfcal{S}^i$ and $\bfcal{P}=\bfcal{S}^e$ obviously result in the lowest cost per iteration (except for $\bfcal{P}=\bfcal{I}$, which corresponds to no preconditioning and is only included for reference), but this does not guarantee \textit{a priori }the lowest total computation time since the number of iterations $R$ required for GMRES convergence also depends on the choice of $\bfcal{P}$, as we demonstrate numerically in \S\ref{sct:benchmarks}.

\begin{table}[t!]
\begin{center}
\begin{tabular}{ll}
\toprule
Choice of $\bfcal{P}$ & Total matvecs\\
\midrule
$\bfcal{I}$ & $(4M^2+4M)(R + \floor*{ R/\rho}) $\\[3pt]
$\bfcal{A}$ & $(8M^2+8M)(R + \floor*{ R/\rho}) + 4M^2+4M$\\[3pt]
$\bfcal{D}$ & $(4M^2+12M)(R + \floor*{ R/\rho}) + 8M$\\[3pt]
$\bfcal{D}^i$, $\bfcal{D}^e$  & $(4M^2+8M)(R + \floor*{ R/\rho}) + 4M$\\[3pt]
$\bfcal{S}^i$, $\bfcal{S}^e$  & $(4M^2+6M)(R + \floor*{ R/\rho}) + 2M$\\[3pt]
\bottomrule
\end{tabular}
\caption{\label{tab:matvecs} Total GMRES solver time in terms of matvecs using different choices of preconditioning operator $\bfcal{P}$. 
Here $R$ is the number of GMRES iterations required to achieve convergence at a specified tolerance (note that $R$ is expected to depend strongly on the choice of $\bfcal{P}$), $\rho$ is the number of iterations per cycle passed as the restart argument in GMRES, $\floor*{\cdot}$ is the ``floor'' function, and $M$ is the number of scatterers. For simplicity, we are assuming that the cost of the boundary integral operators is the same for each scatterer $\Gamma_m$ (or pair of scatterers $(\Gamma_\ell,\Gamma_m)$ in the case of off-diagonal operators such as $\mathcal{C}^e_{\ell,m}$), which will be true provided each scatterer is discretised with roughly the same number of degrees of freedom.
}
\end{center}
\end{table}

\subsection{Discretisation of the operator product}
We now turn to the choice of discretisation strategy for the operator product $\bfcal{P}\bfcal{A}$. 
To set the scene 
we recall from \cite{buffa2003galerkin} that $\mathbf{H}^{-\frac{1}{2}}_\times (\textnormal{div}_{\Gamma_m}, \Gamma_m)$ is self-dual with respect to the anti-symmetric $L^2$ dual pairing\footnote{We recall that saying a Hilbert space $Y$ is dual to another Hilbert space $X$ with respect to a pairing $\langle \cdot, \cdot\rangle_{Y\times X}$ means the map $\mathcal{I}:Y\to X^*$ defined by $\mathcal{I}y (x)=\langle y, x\rangle_{Y\times X}$ is an isomorphism (continuous linear bijection).} 
\begin{align}
\label{dualpairingGammam}
    \langle \mathbf{a},\mathbf{b} \rangle_{\Gamma_m} = \int_{\Gamma_m} \mathbf{a}\cdot (\mathbf{n}_m\times \mathbf{b})\,\mathrm{d}s,
\end{align}
which implies that $\mathbf{X}$ is self-dual with respect to the pairing
\begin{align}
\label{dualdef}
\left\langle \bigoplus_{m=1}^M\left(\begin{array}{c}\mathbf{c}_m \\ \mathbf{d}_m \end{array}\right),\bigoplus_{m=1}^M\left(\begin{array}{c}\mathbf{e}_m \\ \mathbf{f}_m \end{array}\right)\right\rangle = 
\sum_{m=1}^M \langle \mathbf{c}_m ,\mathbf{f}_m \rangle_{\Gamma_m} + \langle\mathbf{d}_m , \mathbf{e}_m \rangle_{\Gamma_m}.
\end{align}
Then for a stable Galerkin discretisation of $\bfcal{P}\bfcal{A}$ we follow the procedure outlined in \cite{betcke2020product}. For each choice of $\bfcal{A}$ and $\bfcal{P}$ we select a discrete domain (trial) space, dual to the range (test) space, and range space, which we denote by ($\mathbf{X}^{\bfcal{A},{\rm dom}}_h$, $\mathbf{X}^{\bfcal{A},{\rm dual}}_h$, $\mathbf{X}^{\bfcal{A},{\rm ran}}_h$) and ($\mathbf{X}^{\bfcal{P},{\rm dom}}_h$, $\mathbf{X}^{\bfcal{P},{\rm dual}}_h$, $\mathbf{X}^{\bfcal{P},{\rm ran}}_h$) respectively. All these spaces should be finite-dimensional subspaces of $\mathbf{X}$ with the same common dimension $N\in \mathbb{N}$, and it should hold that $\mathbf{X}^{\bfcal{A},{\rm dual}}_h$ is dual to $\mathbf{X}^{\bfcal{A},{\rm ran}}_h$ and that $\mathbf{X}^{\bfcal{P},{\rm dual}}_h$ is dual to $\mathbf{X}^{\bfcal{P},{\rm ran}}_h$ with respect to the pairing \eqref{dualdef}. Furthermore, we require that $\mathbf{X}^{\bfcal{A},{\rm ran}}_h = \mathbf{X}^{\bfcal{P},{\rm dom}}_h$. Then the discrete strong form of the operator product $\bfcal{P}\bfcal{A}$ is well defined as a map from $\mathbf{X}^{\bfcal{A},{\rm dom}}_h$ to $\mathbf{X}^{\bfcal{P},{\rm ran}}_h$. 
Choosing bases $\{\bm{\phi}^{\bfcal{A},{\rm dom}}_j\}_{j=1}^N$, $\{\bm{\phi}^{\bfcal{A},{\rm dual}}_j\}_{j=1}^N$, etc.\ for the spaces $\mathbf{X}^{\bfcal{A},{\rm dom}}_h$, $\mathbf{X}^{\bfcal{A},{\rm dual}}_h$, etc., the matrix associated with this operator product is given by $\mathbf{M}_{\bf{P}}^{-1}\mathbf{P}\mathbf{M}_{\bf{A}}^{-1}\mathbf{A}$, where $\mathbf{A}$ and $\mathbf{P}$ are the Galerkin matrices and 
$\mathbf{M}_{\bf{A}}$ and $\mathbf{M}_{\bf{P}}$ are the mass matrices for $\bfcal{A}$ and $\bfcal{P}$ respectively. Explicitly,
\begin{alignat}{3}
\label{GalA}
\mathbf{A}_{ij} &= \langle \bm{\mathcal{A}}\bm{\phi}^{\bfcal{A},{\rm dom}}_j,\bm{\phi}^{\bfcal{A},{\rm dual}}_i \rangle, &&\qquad i,j\in\{1,\ldots,N\},\\
\label{GalP}
\mathbf{P}_{ij} &= \langle \bm{\mathcal{P}}\bm{\phi}^{\bfcal{P},{\rm dom}}_j,\bm{\phi}^{\bfcal{P},{\rm dual}}_i \rangle, &&\qquad i,j\in\{1,\ldots,N\},\\
\label{MassA}
(\mathbf{M}_{\bf{A}})_{ij} &= \langle \bm{\phi}^{\bfcal{A},{\rm ran}}_j,\bm{\phi}^{\bfcal{A},{\rm dual}}_i \rangle,&&\qquad i,j\in\{1,\ldots,N\},\\
\label{MassP}
(\mathbf{M}_{\bf{P}})_{ij} &= \langle \bm{\phi}^{\bfcal{P},{\rm ran}}_j,\bm{\phi}^{\bfcal{P},{\rm dual}}_i \rangle, &&\qquad i,j\in\{1,\ldots,N\}.
\end{alignat}
Given this framework, to obtain an approximate solution of \eqref{eqn:PMCHWT_multiple_scatterers_Precond} we seek $\mathbf{u}^s_h = \sum_{j=1}^N x_j \bm{\phi}^{\bfcal{A},{\rm dom}}_j\in \mathbf{X}^{\bfcal{A},{\rm dom}}_h$ satisfying the linear system
\begin{align}
\mathbf{M}_{\bf{P}}^{-1}\mathbf{P}\mathbf{M}_{\bf{A}}^{-1}\mathbf{A}\mathbf{x} = \mathbf{M}_{\bf{P}}^{-1}\mathbf{P}\mathbf{M}_{\bf{A}}^{-1}\mathbf{b},
\end{align} 
where $\mathbf{x} = (x_1,\ldots,x_N)^T$ and 
$\mathbf{b} = (b_1,\ldots,b_N)^T$, with 
$b_j=\langle \left(\frac{1}{2}\bm{\mathcal{ {I}}} - \bm{\mathcal{{D}}}^i \right) \mathbf{u}^{inc},\bm{\phi}^{\bfcal{A},{\rm dual}}_j \rangle$, $j=1,\ldots,N$.

Next we discuss the selection of the discrete spaces, under the assumption that each $\Gamma_m$ is polyhedral. Introducing a triangulation on $\Gamma_m$, the primal mesh, we define the space RWG$_m$ to be the span of the Rao-Wilton-Glisson basis functions \cite{rao1982electromagnetic} on the primal mesh, and the space BC$_m$ to be the span of the Buffa-Christiansen basis functions \cite{buffa2007dual} on the corresponding dual mesh. The spaces RWG$_m$ and BC$_m$ have the same number of global degrees of freedom, and moreover --- this being the reason for the introduction of the space BC$_m$ in \cite{buffa2007dual} --- BC$_m$ is dual to RWG$_m$ with respect to the pairing \eqref{dualpairingGammam}. This means that the requirements for a stable discretisation of $\bfcal{P}\bfcal{A}$ detailed above can be met by choosing
\begin{align}
    \label{eqn:RWG_discretisation1}\mathbf{X}^{\bfcal{A},{\rm dom}}_h &=\mathbf{X}^{\bfcal{A},{\rm dual}}_h =\mathbf{X}^{\bfcal{P},{\rm ran}}_h =\bigoplus_{m=1}^M \left(\begin{array}{c}{\rm RWG}_m \\ {\rm RWG}_m \end{array}\right), 
    \quad
    \mathbf{X}^{\bfcal{P},{\rm dom}}_h = \mathbf{X}^{\bfcal{P},{\rm dual}}_h =\mathbf{X}^{\bfcal{A},{\rm ran}}_h= \bigoplus_{m=1}^M \left(\begin{array}{c}{\rm BC}_m \\ {\rm BC}_m \end{array}\right), 
\end{align}
which is the choice we make in this paper. We note, however, that other choices are possible, such as 
\begin{align}
    \label{eqn:mixed_discretisation}\mathbf{X}^{\bfcal{A},{\rm dom}}_h &=\mathbf{X}^{\bfcal{A},{\rm dual}}_h =\mathbf{X}^{\bfcal{A},{\rm ran}}_h=\mathbf{X}^{\bfcal{P},{\rm dom}}_h = \mathbf{X}^{\bfcal{P},{\rm dual}}_h  =\mathbf{X}^{\bfcal{P},{\rm ran}}_h  = \bigoplus_{m=1}^M \left(\begin{array}{c}{\rm RWG}_m \\ {\rm BC}_m \end{array}\right), 
\end{align}
which is the choice introduced in \cite{scroggs2017software} and used in \cite{kleanthous2019calderon}. 

Both \eqref{eqn:RWG_discretisation1} and \eqref{eqn:mixed_discretisation} possess symmetries that allow the overall assembly time and memory cost to be reduced by judicious caching of the discretised operators. For \eqref{eqn:mixed_discretisation}, the fact that the same discrete spaces are used for $\bfcal{P}$ and $\bfcal{A}$ means that the assembly of $\mathbf{P}$ comes ``for free'' once $\mathbf{A}$ is assembled, since one can just re-use the appropriate sub-blocks of $\mathbf{A}$ to build $\mathbf{P}$ (all our choices of $\bm{\mathcal{P}}$ require a subset of the operators appearing in $\bm{\mathcal{A}}$). 
For \eqref{eqn:RWG_discretisation1} this is not possible, since different discrete spaces are used for $\bfcal{P}$ and $\bfcal{A}$. However, for \eqref{eqn:RWG_discretisation1}, the fact that the two entries of each factor in the discrete spaces coincide means that double occurrences of identical boundary integral operators in the structure of $\bfcal{A}$ 
at the continuous level (in the diagonal and off-diagonal entries of the matrices in \eqref{AimAem}) give rise to double occurrences of identical discretised versions at the discrete level, which only need to be assembled once, cutting the assembly cost of $\mathbf{A}$ in half compared to a naive implementation. The same is true for $\mathbf{P}$ for the choices of $\bfcal{P}$ presented above. This symmetry does not hold for \eqref{eqn:mixed_discretisation}.

Our reasons for favouring \eqref{eqn:RWG_discretisation1} over \eqref{eqn:mixed_discretisation} in the present study are closely related to the above observations, and to the fact that assembling an operator with BC$_m$ as the trial and/or test space is significantly more expensive than using RWG$_m$, all other factors being equal, because of the need to work with the barycentric mesh, which has six times as many elements as the primal mesh.
In \eqref{eqn:RWG_discretisation1} we are placing this cost burden on $\mathbf{P}$ not $\mathbf{A}$ (as would be the case in \eqref{eqn:mixed_discretisation}). This means that by using a ``reduced'' choice of $\bfcal{P}$ (i.e., $\bfcal{D}$, $\bfcal{D}^e$, $\bfcal{D}^i$, $\bfcal{S}^e$ or $\bfcal{S}^i$) we avoid having to assemble all of the BC$_m$ versions of the operators in $\bfcal{A}$. 
Furthermore, as we discuss in the next section, the use of different discrete spaces for $\bfcal{P}$ and $\bfcal{A}$ in \eqref{eqn:RWG_discretisation1} gives us flexibility to use assembly routines of differing accuracy for $\mathbf{P}$ and $\mathbf{A}$, making it possible to use a lower quality (cheaper) assembly routine for $\mathbf{P}$ than for $\mathbf{A}$, further mitigating the high cost of the BC$_m$ spaces. 

\subsection{Bi-parametric matrix assembly}\label{sct:accelerated_bi_parametric}
To reduce the cost of assembling and storing the matrices $\mathbf{P}$ and $\mathbf{A}$ when the dimension $N$ is large we use hierarchical matrix ($\mathcal{H}$-matrix) compression \cite{hackbusch2015hierarchical}. 
In a recent study \cite{escapil2019fast} of Calder\'on preconditioning for perfectly conducting scatterers using the electric field integral equation (EFIE), it was shown that a significant saving in memory and computation time can be obtained by using a cheaper (lower-quality) $\mathcal{H}$-matrix approximation and quadrature routine for the preconditioner than for the original operator, termed a ``bi-parametric'' approach in \cite{escapil2019fast}. 
Its theoretical foundations can be found in \cite{ESCAPILINCHAUSPE2021220}.
We aim to show in the current paper that the same is true for the PMCHWT formulation for (multiple) dielectric scatterers. 

Before specifying our bi-parametric approach we describe some relevant features of the $\mathcal{H}$-matrix implementation and quadrature routines in the software package Bempp with which we obtain our numerical results in \S\ref{sct:benchmarks} and \S\ref{sec:Applications}. Since $\mathcal{H}$-matrix implementations and BEM quadrature routines are complicated and vary considerably, we do not provide full details of our algorithms, but rather highlight the key general principles, with the expectation that expert readers will adapt our ideas to their own specific implementations.

Given index sets $\mathcal{I}= \mathcal{I}^{(0)}_1=\{1, \ldots, N\}$ and $\mathcal{J}= \mathcal{J}^{(0)}_1=\{1, \ldots, N\}$ for the degrees of freedom (dofs) in the test and trial spaces respectively, we first generate ``cluster trees'' $\mathcal{T}(\mathcal{I})$ and $\mathcal{T}(\mathcal{J})$ by repeated subdivision of $\mathcal{I}$ and $\mathcal{J}$ according to certain geometric criteria. 
We then generate a ``block cluster tree'' by applying specified geometric admissibility conditions to pairs of nodes from $\mathcal{T}(\mathcal{I})$ and $\mathcal{T}(\mathcal{J})$, each of which corresponds (after suitable re-ordering of the dofs) to a sub-block of the Galerkin matrix. 
In our implementation we adopt the following admissibility condition: a pair of nodes $\mathcal{I}_p^{(q)}$ and $\mathcal{J}_r^{(s)}$ of level $q$ (resp.\ $s$) and index $p$ (resp.\ $r$) are deemed ``admissible'' if certain geometric bounding boxes $X(\mathcal{I}_p^{(q)})$ and $Y(\mathcal{J}_r^{(s)})$ associated with the collections of mesh elements on which the basis functions indexed by $\mathcal{I}_p^{(q)}$ and $\mathcal{J}_r^{(s)}$ are supported satisfy
\begin{align}
    \text{dist}(X(\mathcal{I}_p^{(q)}),Y(\mathcal{J}_r^{(s)})) >0.
\end{align}
If this condition fails, the block is further refined until all its sub-blocks are deemed admissible or are smaller than some specified minimum block size, in which case they are deemed ``inadmissible''. Inadmissible blocks are fully (densely) assembled, while admissible blocks are represented by low-rank approximations, generated using adaptive cross approximation (ACA) \cite{bebendorf2000approximation, kurz2002adaptive, bebendorf2003adaptive}. Given an admissible block with underlying matrix $\mathbf{B}$ and a user-specified tolerance parameter $\nu >0$, the ACA algorithm delivers an approximation
\begin{align}
    \mathbf{B}_\nu = \sum_{i=1}^{r} \mathbf{u}_i \mathbf{v}_i^H,
\end{align}
of rank $r=r(\nu,\mathbf{B})$ such that 
\begin{align}
    \frac{\|\mathbf{B}_\nu - \mathbf{B}\|_2}{\|\mathbf{B}\|_F}  \leq \nu.
\end{align}
Varying $\nu$ allows the user to control the quality (and cost) of the low-rank approximations. 
Our $\mathcal{H}$-matrix implementation has an additional feature, which saves time and memory by assembling only a subset of the admissible blocks corresponding to ``near-field'' interactions. Given a user-specified cutoff parameter $\chinf\in [0,\infty]$, the admissible blocks for which 
\begin{equation}
\label{eqn:nf_condition}\text{dist}(X(\mathcal{I}_p^{(q)}),Y(\mathcal{J}_r^{(s)})) \leq \chinf
\end{equation}
are assembled using ACA, while all other admissible blocks are set to zero. Setting $\chinf = \infty$ corresponds to assembling all admissible blocks, as in a standard implementation. 

Assembly of both the admissible and inadmissible blocks for $\mathbf{A}$ and $\mathbf{P}$ requires the evaluation of the Galerkin integrals \eqref{GalA} and \eqref{GalP} respectively. Each of these can be written in terms of integrals of the form 
\begin{align}
\int_{T_1} \int_{T_2} F(\mathbf{x},\mathbf{y}) \mathrm{d} s(\mathbf{y}) \mathrm{d} s (\mathbf{x}), 
\end{align}
where $T_1$ and $T_2$ are triangles in either the primal mesh or the barycentric mesh, and the integrand $F(\mathbf{x},\mathbf{y})$ involves the fundamental solution $G (\mathbf{x}, \mathbf{y}) = \frac{\exp (\mathrm{i} k | \mathbf{x} - \mathbf{y}|)}{4 \pi | \mathbf{x} - \mathbf{y}|}$ and a pair of discrete basis functions. When the closures of $T_1$ and $T_2$ are disjoint, $F(\mathbf{x},\mathbf{y})$ is smooth and a standard tensor product Gauss rule based on symmetric Gauss points over triangles is used. The number of quadrature points used depends on whether the integral is classified by Bempp as either near-, medium- or far field (based on the distance between the triangles and their sizes), with the user specifying quadrature order parameters $q_{near}$, $q_{medium}$ and $q_{far}$ for each case. 
When $T_1$ and $T_2$ share a vertex, edge or are the same triangle, $F(\mathbf{x},\mathbf{y})$ is singular and quadrature routines of the type described in \cite{sauter2010boundary} are applied, with the user specifying a single quadrature order parameter $q_{singular}$ for all singular cases. 

The idea of a so-called ``bi-parametric'' implementation is to use different choices of the parameters $\nu$, $\chinf$ and $\mathbf{q}=(q_{near}, q_{medium}, q_{far}, q_{singular})$ 
for the assembly of $\mathbf{P}$ and $\mathbf{A}$. Naturally we label these $\nu_{\mathbf{P}}$, $\chinf_{\mathbf{P}}$, $\mathbf{q}_{\mathbf{P}}$ and $\nu_{\mathbf{A}}$, $\chinf_{\mathbf{A}}$, $\mathbf{q}_{\mathbf{A}}$, respectively. Our expectation is that accuracy of matrix assembly should be more important for $\mathbf{A}$ than for $\mathbf{P}$, since $\mathbf{A}$ governs overall solution accuracy whereas $\mathbf{P}$ is included merely to accelerate the convergence of the iterative solver. Hence we expect to be able to use a larger value of $\nu_{\mathbf{P}}$, and smaller values of $\chinf_{\mathbf{P}}$ and $\mathbf{q}_{\mathbf{P}}$, compared to $\nu_{\mathbf{A}}$, $\chinf_{\mathbf{A}}$ and $\mathbf{q}_{\mathbf{A}}$, saving computation time and memory. We demonstrate that this is possible in a set of numerical experiments in \S\ref{sct:benchmarks} and \S\ref{sec:Applications}.

\begin{figure}[t!]
    \centering
    \includegraphics[width = 0.24 \textwidth]{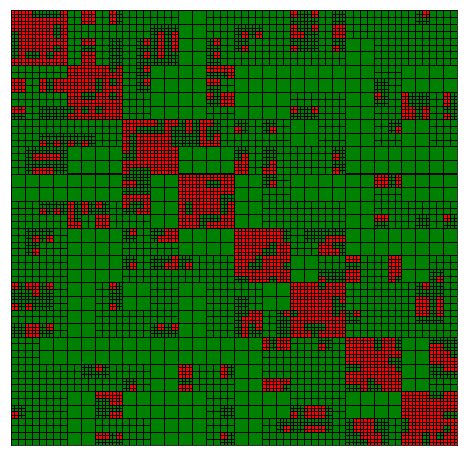}
    \includegraphics[width = 0.24 \textwidth]{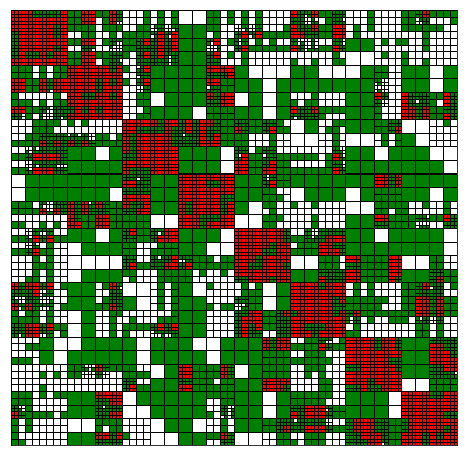}
    \includegraphics[width = 0.24 \textwidth]{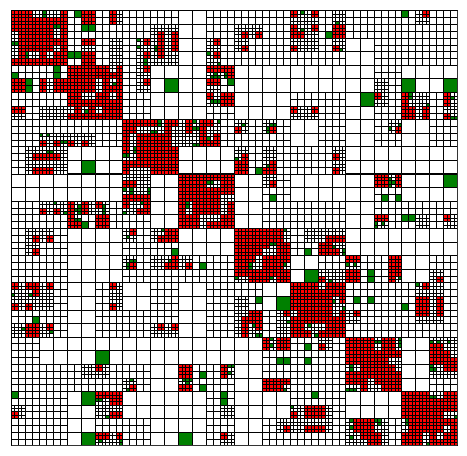}
    \includegraphics[width = 0.24 \textwidth]{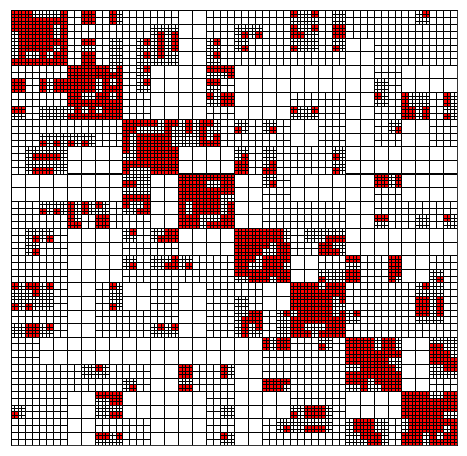}\\
 $\chinf=\infty$ 
\hspace{26mm} $\chinf=0.5$
\hspace{26mm} $\chinf=0.1$
\hspace{26mm} $\chinf=0$
    \caption{Block cluster trees produced by Bempp in the assembly of a single electric field operator $\mathcal{S}$ on the unit cube $[0,1]^3$ with $k=5$ and approximately 10 elements per wavelength in the BEM mesh.  Here red indicates inadmissible blocks, green indicates admissible blocks that require ACA approximation, and white indicates admissible blocks that do not require assembly. As the near-field cutoff parameter $\chinf$ decreases from left to right the overall compression rate decreases, taking the values 0.83, 0.60, 0.17 and 0.14 respectively when the ACA parameter $\nu=0.001$.}      
    \label{fig:H_matrix_nf}
\end{figure}

To give an idea of the potential savings of this bi-parametric approach, an illustration of the effect of changing the near-field cutoff parameter $\chi$ on the cluster tree is shown in Figure \ref{fig:H_matrix_nf}. 
As expected, the overall compression ratio (the memory cost as a fraction of that of dense assembly) drops as we remove more and more of the far-field. We chose to present a relatively low-frequency case in Figure \ref{fig:H_matrix_nf} to limit the number of degrees of freedom so as to make the block structure easier to see, but for the high frequency applications in \S\ref{sec:Applications} the compression ratios are significantly smaller. 
Regarding quadrature, 
for non-singular integrals the default values $(q_{near}, q_{medium}, q_{far})=(4,3,2)$ correspond to $(36,16,9)$ integrand evaluations respectively. Reducing these values to the minimum possible $(q_{near}, q_{medium}, q_{far})=(1,1,1)$ brings the numbers of integrand evaluations down to $(1,1,1)$ respectively. 
For singular integrals, the default value $q_{singular}=6$ means 512 integrand evaluations if $T_1$ and $T_2$ share a single common vertex, 1280 if $T_1$ and $T_2$ share a single common edge, and 1536 if $T_1$ and $T_2$ coincide. Reducing this value to the minimal possible $q_{singular}=1$ brings these totals down to 2, 5 and 6 respectively. 

\section{Numerical results}\label{sct:benchmarks}
\begin{figure}[tp!]
    \centering
    \includegraphics[width = 0.5\textwidth]{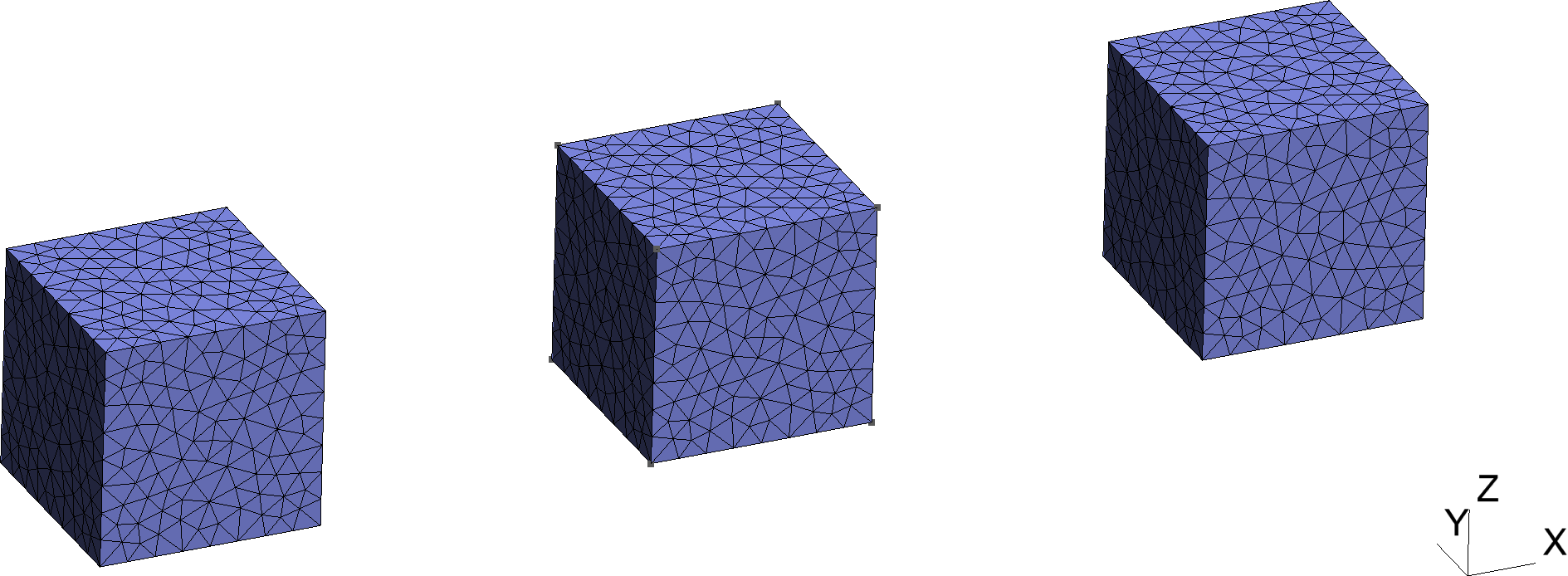} \hfill    
    \includegraphics[width = 0.45\textwidth]{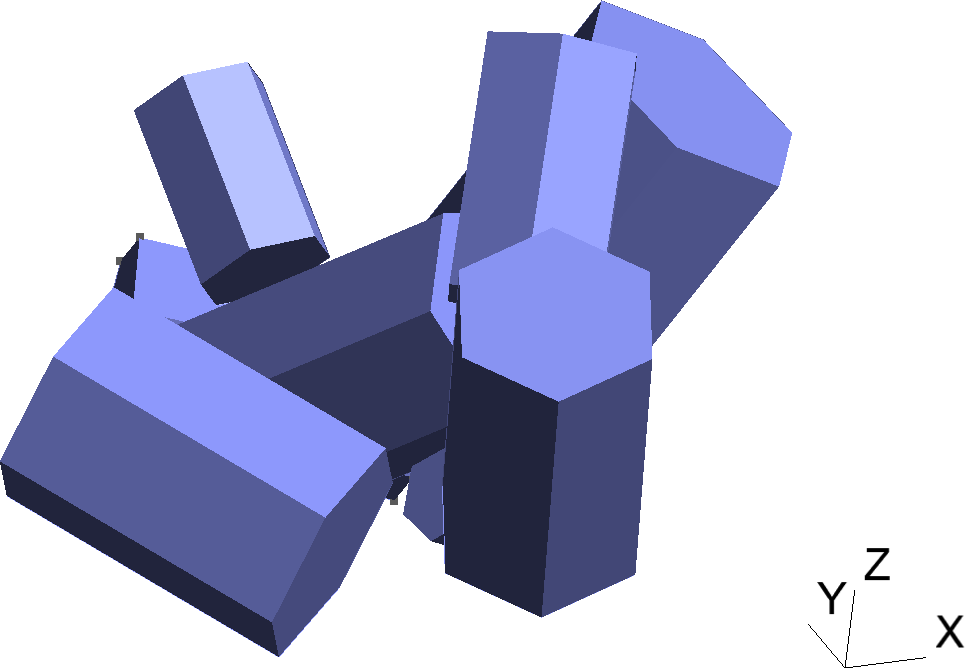}
    \caption{
    \label{fig:geometry}  Scatterer configurations. Left: 3 cubes. Right: 8-branch ice crystal aggregate from \cite{yang1998single}.}
\end{figure}

In this section we investigate the performance of the acceleration techniques described in \S\ref{sct:calderon_construction}
in the context of a simple benchmark problem, namely scattering of an incident plane wave $\mathbf{E}^{inc} = (0,0, \mathrm{e}^{\mathrm{i}k_e x})$ by an array of 3 disjoint identical cubes, as illustrated in Figure \ref{fig:geometry} (left panel). The cubes have side length $0.4$, are aligned with the coordinate axes, and have their front bottom left vertices at the points $(-1,0,0)$, $(0,0,0)$ and $(1,0,0)$ respectively. We take the refractive index of each cube to be $n = 1.311 + 2.289 \times 10^{-9}\mathrm{i}$, which is a representative value for our intended atmospheric physics application, being the measured value of the refractive index of ice at wavelength $\lambda = 0.55$ $\mu$m \cite{warren2008optical}. 

In all our experiments in this section and the next, scatterers are meshed using a maximum mesh size of $h = 2 \pi/(10 k_e)$, ensuring at least 10 elements per wavelength in each coordinate direction. For the discretisation of the original operator $\bfcal{A}$ we always use $\nu_{\mathbf{A}}=0.001$, $\chinf_{\mathbf{A}}=\infty$, $\mathbf{q}_{\mathbf{A}}=(4,3,2,6)$, which are the default Bempp parameter values. Linear systems are solved with restarted GMRES with parameters $tol = 10^{-5}$ and $restart = 200$.
A large $restart$ parameter was chosen to remove the influence of restarting on convergence. Hence, any observations on the improvement of convergence are due to the effectiveness of the reduced preconditioners, bi-parametric assembly, or combination of the two rather than the GMRES parameters.
For the problems considered here, comparison to a reference solution obtained via dense assembly (where possible) suggests that this produces solutions accurate to within about 0.2\%.

\subsection{Reduced Calder\'on preconditioning}
We begin by comparing the performance of the block-diagonal preconditioner $\bfcal{P}=\bm{\mathcal{D}}$ to that of the various ``reduced'' preconditioners $\bfcal{P}=\bm{\mathcal{D}}^e$, $\bfcal{P}=\bm{\mathcal{D}}^i$, $\bfcal{P}=\bm{\mathcal{S}}^e$ and $\bfcal{P}=\bm{\mathcal{S}}^i$ introduced in \S\ref{sec:Pchoice}. Throughout this subsection we use the default parameters $\nu_{\mathbf{P}}=0.001$, $\chinf_{\mathbf{P}}=\infty$, $\mathbf{q}_{\mathbf{P}}=(4,3,2,6)$. Bi-parametric implementations are not considered until the next subsection.
For completeness, we include information about the performance of the non preconditioned operator $\bfcal{A}$ (without mass-matrix preconditioning) where relevant.

\begin{figure}[th!]
    \centering
    \includegraphics[width = \textwidth]{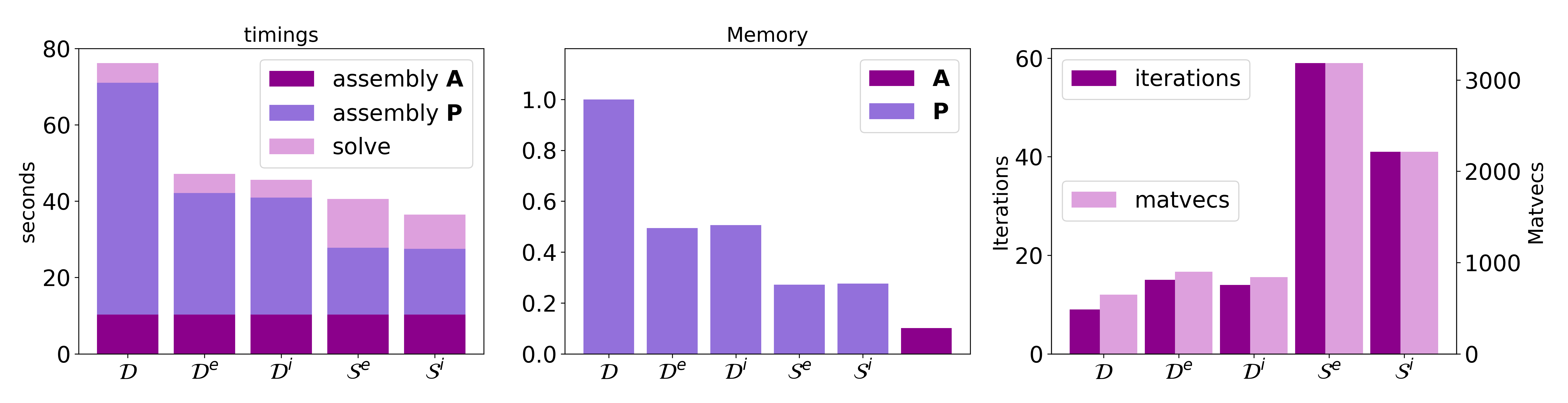} \\
    \includegraphics[width = 0.8 \textwidth]{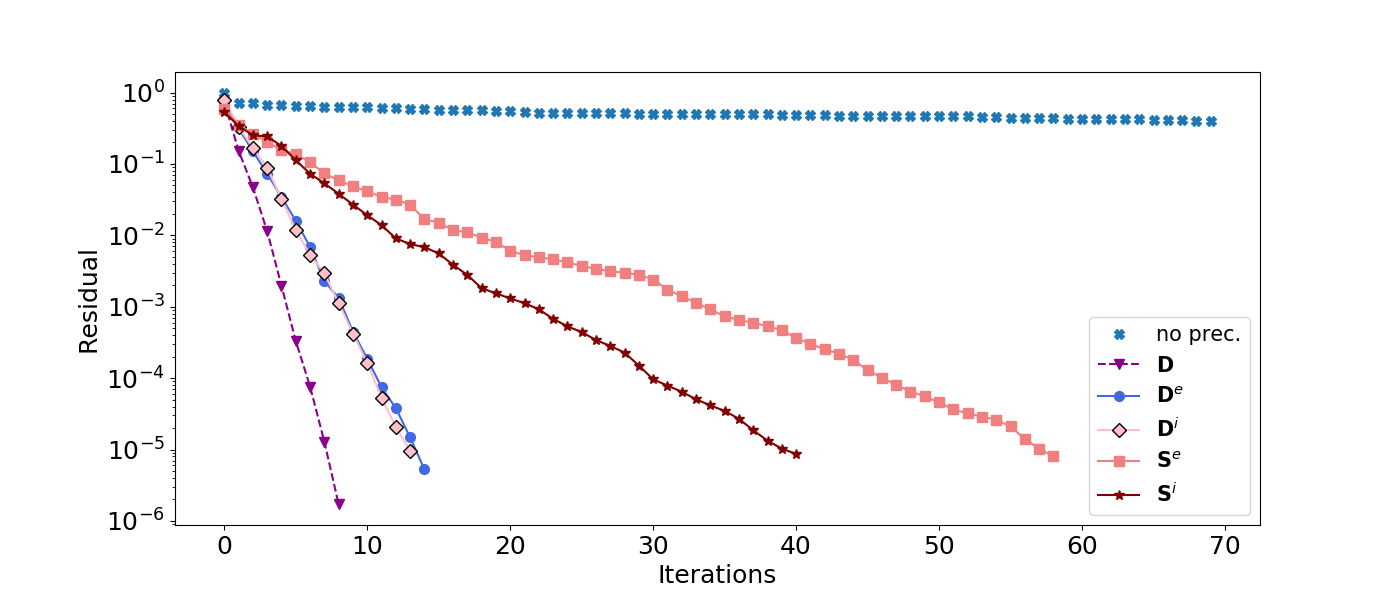}
    \caption{Computation times, memory costs (normalised relative to the cost of $\bfcal{D}$),  GMRES iteration/matvec counts and residual curves for the reduced preconditioners $\bfcal{D}^e$, $\bfcal{D}^i$, $\bfcal{S}^e$ and $\bfcal{S}^i$ (and for the non-preconditioned case in the bottom figure), for scattering by the three cubes in Figure \ref{fig:geometry} (left panel) with $k_e = 11.4$.
    We note that the GMRES iteration for the non-preconditioned operator was terminated after 2000 iterations (without having achieved the target tolerance), at which point it had already exceeded the total time taken for assembly and solution for $\bfcal{DA}$.}
    \label{fig:reducedCalderon_single_k}
\end{figure}

In Figure \ref{fig:reducedCalderon_single_k}, we report assembly and solve times, memory consumption and GMRES iteration and matvec counts for each of the reduced preconditioners in the case $k_e = 11.4$, corresponding to an exterior wavelength $\lambda_e = 2\pi/k_e=0.55$ (making the sides of each cube about 0.7 exterior wavelengths long) and a total number of dofs $N=4395$. As predicted in \S\ref{sec:Pchoice}, the assembly time and memory cost of $\bfcal{S}^i$ and $\bfcal{S}^e$ are half those of $\bfcal{D}^i$ and $\bfcal{D}^e$, which in turn are half those of $\bfcal{D}$. The solve time for $\bfcal{S}^i$ and $\bfcal{S}^e$ is longer than for $\bfcal{D}^i$, $\bfcal{D}^e$, and $\bfcal{D}$, due to a significantly increased iteration and matvec count (i.e.\ weakened preconditioning effect). However, in terms of total computation time $\bfcal{S}^i$ performs best at this value of $k_e$, with $\bfcal{S}^e$, $\bfcal{D}^i$ and $\bfcal{D}^e$ close behind, because the high computational cost of the barycentric refinement means that assembly time for $\bfcal{P}$ significantly dominates solver time in this instance. 

\begin{figure}[th!]
    \centering
    \includegraphics[width = \textwidth]{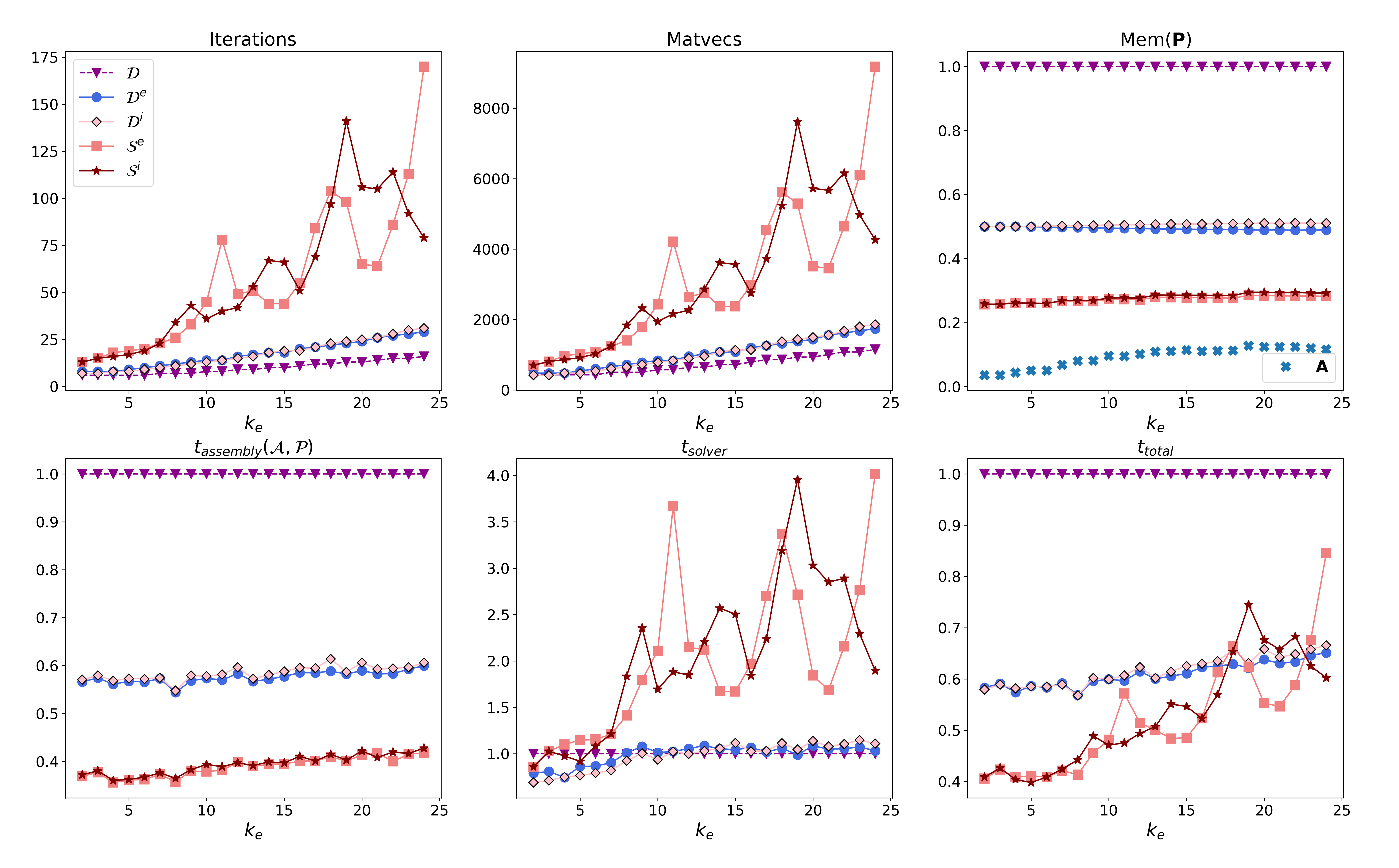}
    \caption{Performance of the reduced preconditioners $\bfcal{D}^e$, $\bfcal{D}^i$, $\bfcal{S}^e$ and $\bfcal{S}^i$ as a function of the exterior wavenumber $k_e$, for scattering by the three cubes in Figure \ref{fig:geometry} (left panel). Memory costs are for the preconditioner and are normalised relative to that of $\bfcal{D}$. We also include the memory cost of operator $\mathbf{A}$ as a reference. Assembly times correspond to the time taken to assemble both $\bfcal{P}$ and $\bfcal{A}$. All timings are normalised relative to those for $\bfcal{D}$. Note that absolute memory costs and timings (not reported here) grow with increasing $k_e$ because the number of degrees of freedom grows approximately quadratically with $k_e$ in order to maintain discretisation quality. As a reference, the non-preconditioned operator $\bfcal{A}$ required $662, 484, 758, 1734$, and $1682$ GMRES iterations for the first five wavenumbers, with normalised $t_{solver}$ at 11.4, 8.6, 10.3, 24.9, 26.7, and $t_{total}$ 1.0,  0.7,  0.8, 1.7, 1.7 respectively. For the remaining wavenumbers GMRES was stopped after 2000 iterations without converging to the desired tolerance. }
    \label{fig:reducedCalderon_multiple_vs_k}
\end{figure}

In Figure \ref{fig:reducedCalderon_multiple_vs_k}, we present similar data for the range of wavenumbers $k_e=1,2,\ldots 24$. In our simulations, absolute memory costs and timings (not reported here) grow with increasing $k_e$ because our prescription that $h = 2 \pi/(10 k_e)$ means that the number of degrees of freedom grows approximately quadratically with $k_e$. However, to allow easier comparison between the different preconditioners, for each fixed $k_e$ we have normalised memory costs and timings relative to those for $\bfcal{P}=\bfcal{D}$. From the third panel we see that, as expected, the memory cost for $\bfcal{D}^i$ and $\bfcal{D}^e$ is half that of $\bfcal{D}$, and that of $\bfcal{S}^i$ and $\bfcal{S}^e$ is half that again. In the fourth panel we see a similar trend in assembly time, although the improvement factors are slightly worse than a half in each case because we are reporting total assembly times for $\bfcal{P}$ and $\bfcal{A}$, not just for $\bfcal{P}$ alone as in Figure \ref{fig:reducedCalderon_single_k}. The number of GMRES iterations (and matvecs) required for $\bfcal{D}^i$ and $\bfcal{D}^e$ is roughly double that for $\bfcal{D}$, while for $\bfcal{S}^i$ and $\bfcal{S}^e$ the iteration and matvec count is much higher, and somewhat erratic. 
This is reflected in the solver time, which for $\bfcal{D}^i$ and $\bfcal{D}^e$ is lower than that of $\bfcal{D}$ for small $k_e$ and essentially the same as that of $\bfcal{D}$ for $k_e$ greater than $10$, and for $\bfcal{S}^i$ and $\bfcal{S}^e$ is up to four times than of $\bfcal{D}$ in the range of $k_e$ studied. The total time (assembly plus solver) for $\bfcal{D}^i$ and $\bfcal{D}^e$ is consistently around 60\% that of $\bfcal{D}$, while for $\bfcal{S}^i$ and $\bfcal{S}^e$ it is around 40\% of that of $\bfcal{D}$ for small $k_e$ but rises above that of $\bfcal{D}^i$ and $\bfcal{D}^e$ for some larger values of $k_e$. 
The erratic behaviour of $\bfcal{S}^i$ and $\bfcal{S}^e$ is associated with the non-invertibility of $\mathcal{S}^i_m$ and $\mathcal{S}^e_m$ (and hence $\bfcal{S}^i$ and $\bfcal{S}^e$) at certain (real) resonant wavenumber \cite{buffa2003galerkin}. 
In our case the preconditioners are invertible because we consider a complex wavenumber (the refractive index is $n = 1.311 + 2.289 \times  10^{-9}\mathrm{i}$), but the imaginary part of the wavenumber is small, so the resonance effects are still noticeable. This could be mitigated by using a larger imaginary part in the wavenumber for $\bfcal{S}^i$ and $\bfcal{S}^e$ (as in, e.g., \cite{contopanagos2002well, antoine2008integral}). 
However, we do not pursue this idea here since, as we show in the next section, the resonance effects are controlled effectively by our bi-parametric implementation.

\subsection{Bi-parametric implementation}
We now investigate the performance of the bi-parametric approach outlined in \S\ref{sct:accelerated_bi_parametric}, in which we use a cheaper (lower-quality) $\mathcal{H}$-matrix assembly routine for $\bfcal{P}$ than for $\bfcal{A}$. We first present results for $\bfcal{P}=\bfcal{D}$, and then for the reduced preconditioners of the previous section, $\bfcal{D}^i$, $\bfcal{D}^e$, $\bfcal{S}^i$, $\bfcal{S}^e$.  

\begin{table}[t!]
\centering
\resizebox{\textwidth}{!}{
\begin{tabular}{|r|r|r|rr|rr|rr|rr|rr|}
\toprule
$\chinf_{\mathbf{P}}$ & $\mu_{\mathbf{P}}$ &  $\mathbf{q}_{\mathbf{P}}$   & \multicolumn{2}{c|}{Iters (Matvecs)}  & \multicolumn{2}{c|}{Mem$(\mathbf{P})$} & \multicolumn{2}{c|}{$t_{assembly}$} & \multicolumn{2}{c|}{$t_{solver}$} & \multicolumn{2}{c|}{$t_{total}$}   \\   
\midrule
$\infty$ & $0.001$ & $(4,3,2,6)$ & \textbf{6} (\textbf{456}) & \textbf{9} (\textbf{672}) & 1 & 1 & 1 & 1 & 1 & 1 & 1 & 1  \\[5pt] 
$\infty$ &0.001 & $(1, 1, 1, 1)$ &   7 (528) & \textbf{9} (\textbf{672}) & 1 & 1 &  0.40 & 0.52 & 1.16 & 1.02 & 0.45 & 0.55\\[5pt]
$\infty$ & 0.01 & $(4, 3, 2, 6)$ & \textbf{6} (\textbf{456}) & \textbf{9} (\textbf{672}) & 0.78 & 0.70 & 0.89 & 0.86 & \textbf{0.99} & 0.92 & 0.90 & 0.87 \\[5pt]
$\infty$ & 0.01 & $(1, 1, 1, 1)$ & 7 (582) & \textbf{9} (\textbf{672}) & 0.78 & 0.71 & 0.37 & 0.41 & 1.15 & 0.92 & 0.43 & 0.44 \\[5pt]
$\infty$ & 0.1 & $(1, 1, 1, 1)$ & 7 (582) & \textbf{9} (\textbf{672}) & 0.57 & 0.44 & 0.36 & 0.40 & 1.14 & \textbf{0.86} &  0.42 & 0.43 \\[5pt]
$\infty$ & 0.5 & $(1, 1, 1, 1)$ & 9 (672) & 11 (816) & 0.41 & 0.31 & 0.33 & 0.35 & 1.44 & 1.02 &  0.41 & \textbf{0.39} \\[5pt]
1.0 & 0.1 & $(1, 1, 1, 1)$ & 7 (582) & \textbf{9} (\textbf{672}) & 0.57 & 0.44 & 0.35 & 0.36 & 1.12 & \textbf{0.86} &  \textbf{0.40} & \textbf{0.39} \\[5pt]
0.1 & 0.1 & $(1, 1, 1, 1)$ & 84 (6072) & 16 (1176) & 0.41 & 0.29 &  0.30 & 0.34 &  13.31 & 1.48 &  1.27 & 0.41 \\[5pt]
0 & 0.1 & $(1, 1, 1, 1)$ & 198 (14280) & 36 (2616) & \textbf{0.36} & \textbf{0.18} & \textbf{0.29} & \textbf{0.29} & 31.16 & 3.18 & 2.58 
& 0.48 \\[5pt]
\midrule
\multicolumn{11}{l}{Results obtained using the mixed discretisation \eqref{eqn:mixed_discretisation} from \cite{kleanthous2019calderon}:}\\[1mm]
\midrule
 $\infty$  &  $0.001$ & $(4, 3, 2, 6)$ & \textbf{5} (\textbf{384}) & \textbf{9} (\textbf{672}) & 2 & 2 & 2.08 & 1.98 & 1.90 & 2.25 & 2.07 & 2.00 \\[5pt] 
\bottomrule
\end{tabular}
}
\caption{Performance of various bi-parametric implementations of the block-diagonal preconditioner $\bfcal{P}=\bfcal{D}$ for scattering by the three cubes in Figure \ref{fig:geometry} (left panel). Results for $k_e = 2.1$ appear in the left sub-columns, followed by results for $k_e =11.4$ in the right sub-columns. Memory cost is for $\mathbf{P}$ alone - the relative cost of $\mathbf{A}$ is 0.02 and 0.10 for the two wavenumbers respectively. Assembly time is for both $\bfcal{P}$ and $\bfcal{A}$ combined. For each wavenumber the memory costs and timings have been normalised relative to those for the first row, which corresponds to a non-bi-parametric implementation with 
$(\nu_{\mathbf{P}},\chinf_{\mathbf{P}},\mathbf{q}_{\mathbf{P}})=(\nu_{\mathbf{A}},\chinf_{\mathbf{A}},\mathbf{q}_{\mathbf{A}})=(0.001,\infty,(4,3,2,6))$.
Minimal values in each column are indicated in bold type. 
For completeness, we also include (in the bottom row) corresponding results for the mixed discretisation \eqref{eqn:mixed_discretisation} used in \cite{kleanthous2019calderon}.}
\label{table:multiple_particle}
\end{table}

In Table \ref{table:multiple_particle}, we report results for the same refractive index but two different exterior wavenumbers, $k_e = 2.1$ and $k_e = 11.4$, corresponding to approximately 0.1 and 0.7 wavelengths along each cube side, and $N=378$ and $N=4395$ dofs respectively. For each wavenumber the memory costs and timings have been normalised relative to those for the first row, which corresponds to a non-bi-parametric implementation with $(\nu_{\mathbf{P}},\chinf_{\mathbf{P}},\mathbf{q}_{\mathbf{P}})=(\nu_{\mathbf{A}},\chinf_{\mathbf{A}},\mathbf{q}_{\mathbf{A}})=(0.001,\infty,(4,3,2,6))$. 
The results in the first six rows of the table show that, with $\chi_{\mathbf{P}}=\infty$ fixed, reducing the quadrature orders to the minimum possible $\mathbf{q}_{\mathbf{P}}=(1,1,1,1)$, and increasing the ACA parameter $\nu_{\mathbf{P}}$ from 0.001 through 0.01 and 0.1 and even to 0.5 gives a significant reduction in memory cost and assembly time, with little or no effect on iteration/matvec count and solver time. Furthermore, the results in these first six rows are similar for the two choices of $k_e$. In rows 7-9 of the table we fix $\nu_{\mathbf{P}}=0.1$ and $\mathbf{q}_{\mathbf{P}}=(1,1,1,1)$ and reduce $\chi_{\mathbf{P}}$, which corresponds to neglecting more and more of the far-field behaviour in the preconditioner. With $\chi_{\mathbf{P}}=1.0$ we see similar behaviour to that with $\chi_{\mathbf{P}}=\infty$. But as $\chi_{\mathbf{P}}$ is reduced to 0.1 and then 0, we see a noticeable reduction in memory cost and assembly time. Unfortunately, this comes at the cost of a significant increase in iteration/matvec count and hence solver time. However, the fact that this increase appears to be more serious for $k_e=2.1$ than for $k_e=11.4$ suggests that neglecting far-field behaviour may be possible provided the frequency is not too low.

For completeness, in the final row of the Table \ref{table:multiple_particle} we report performance statistics for $\bfcal{P}=\bfcal{D}$ using the mixed discretisation \eqref{eqn:mixed_discretisation} of \cite{kleanthous2019calderon}. For both wavenumbers considered, while this method performs well in terms of iteration/matvec count, in terms of computation time and memory cost it is roughly twice as expensive as our reference method (i.e. the top row of Table \ref{table:multiple_particle}), because of the high cost of the barycentric refinement. 
For this method we use the default parameters, since a bi-parametric approach is not appropriate, given that the operators in the preconditioner can be simply re-used from the original operator (as explained in \S\ref{sct:calderon_construction}). 
For this same reason, one can argue that the memory costs reported for this method do not provide a fair comparison with the reference method, since storing the preconditioner is actually ``free' once the operator has been assembled and cached. However, when one performs the arguably fairer comparison of dividing the cost of storing the operator with the mixed discretisation \eqref{eqn:mixed_discretisation} by the total cost of storing both the preconditioner \textit{and }the operator with our default discretisation \eqref{eqn:RWG_discretisation1}, one obtains the ratios $2.09$ and $1.88$, for the two wavenumbers respectively. So even in this measure the memory cost is roughly twice that of our reference method.

\begin{figure}[thp!]
    \centering
    \includegraphics[width = \textwidth]{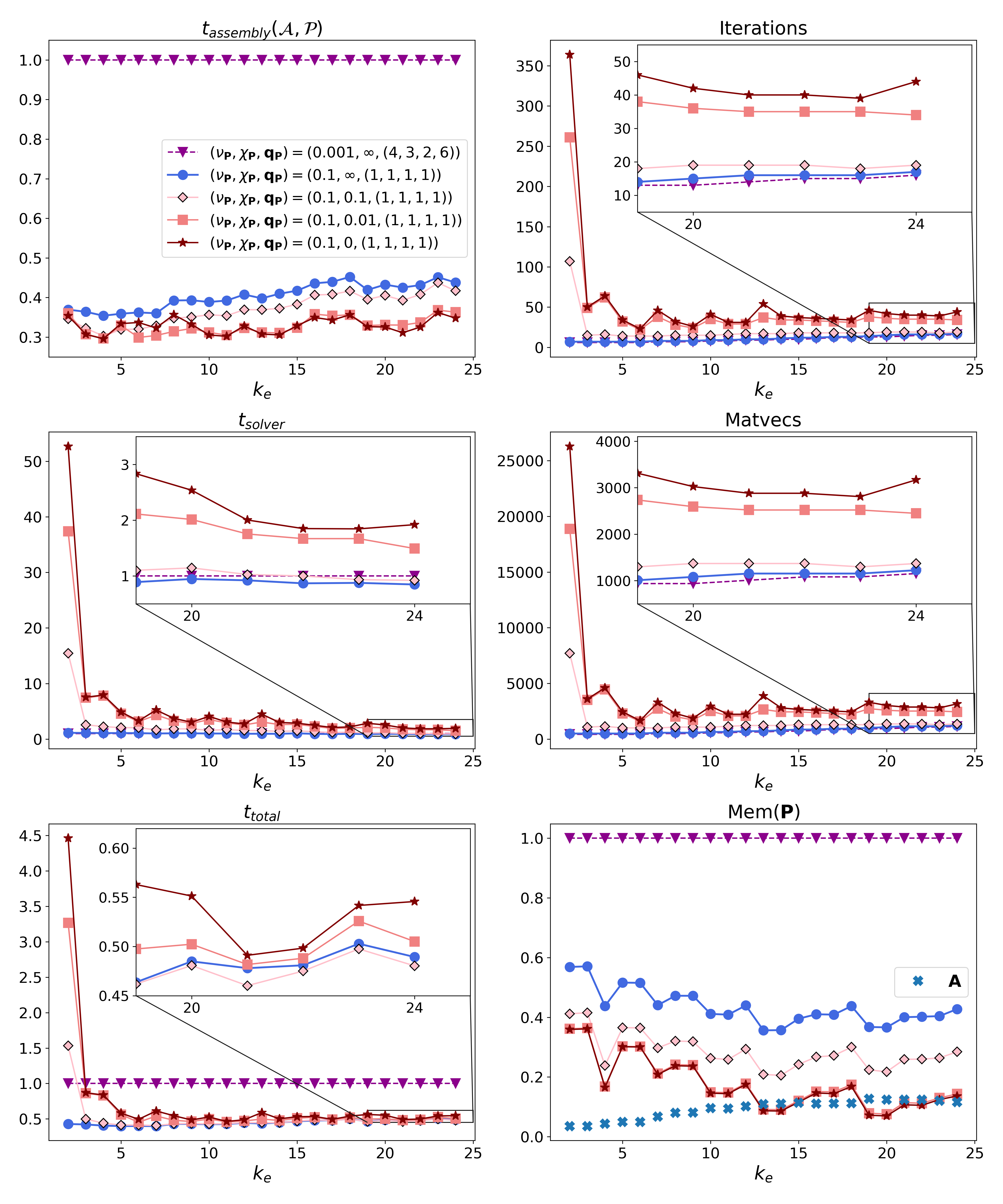}
    \caption{Performance of bi-parametric implementations of the block-diagonal preconditioner $\bfcal{P}=\bfcal{D}$ as a function of $k_e$, for scattering by the three cubes in Figure \ref{fig:geometry} (left panel). The focus is on the effect of varying the near-field cutoff parameter $\chi_{\mathbf{P}}$. Memory costs are for $\mathbf{P}$ alone, while assembly times are for both $\bfcal{P}$ and $\bfcal{A}$ combined. We also include the memory cost of operator $\mathbf{A}$ as a reference. For each wavenumber, memory costs and timings have been normalised relative to those for a non-bi-parametric implementation with $\nu_{\mathbf{P}}=\nu_{\mathbf{A}}=0.001$, $\chinf_{\mathbf{P}}=\chinf_{\mathbf{A}}=\infty$, $\mathbf{q}_{\mathbf{P}}=\mathbf{q}_{\mathbf{A}}=(4,3,2,6)$.} \label{fig:biparametric_D_vs_k}
\end{figure}

To investigate this apparent frequency dependence further, in Figure \ref{fig:biparametric_D_vs_k} we report results for $\bfcal{P}=\bfcal{D}$ for the range of wavenumbers $k_e=1,2,\ldots 24$, with $\nu_{\mathbf{P}}=0.1$ and $\mathbf{q}_{\mathbf{P}}=(1,1,1,1)$ fixed and $\chi_{\mathbf{P}}=\infty$, $1$, $0.1$, $0.01$, and $0$. These results confirm that the significant increase in iteration/matvec count for small values of $\chi_{\mathbf{P}}$ is a low-frequency issue, and that at larger values of $k_e$ the increase in solver time is modest, even for $\chi_{\mathbf{P}}=0$, and is balanced out by the decrease in assembly time, to give a total time only slightly larger than for $\chi_{\mathbf{P}}=\infty$, but with approximately half the memory cost (approximately 20\% of that for the non-bi-parametric reference case). 

\begin{figure}[thp!]
    \centering
    \includegraphics[width = \textwidth]{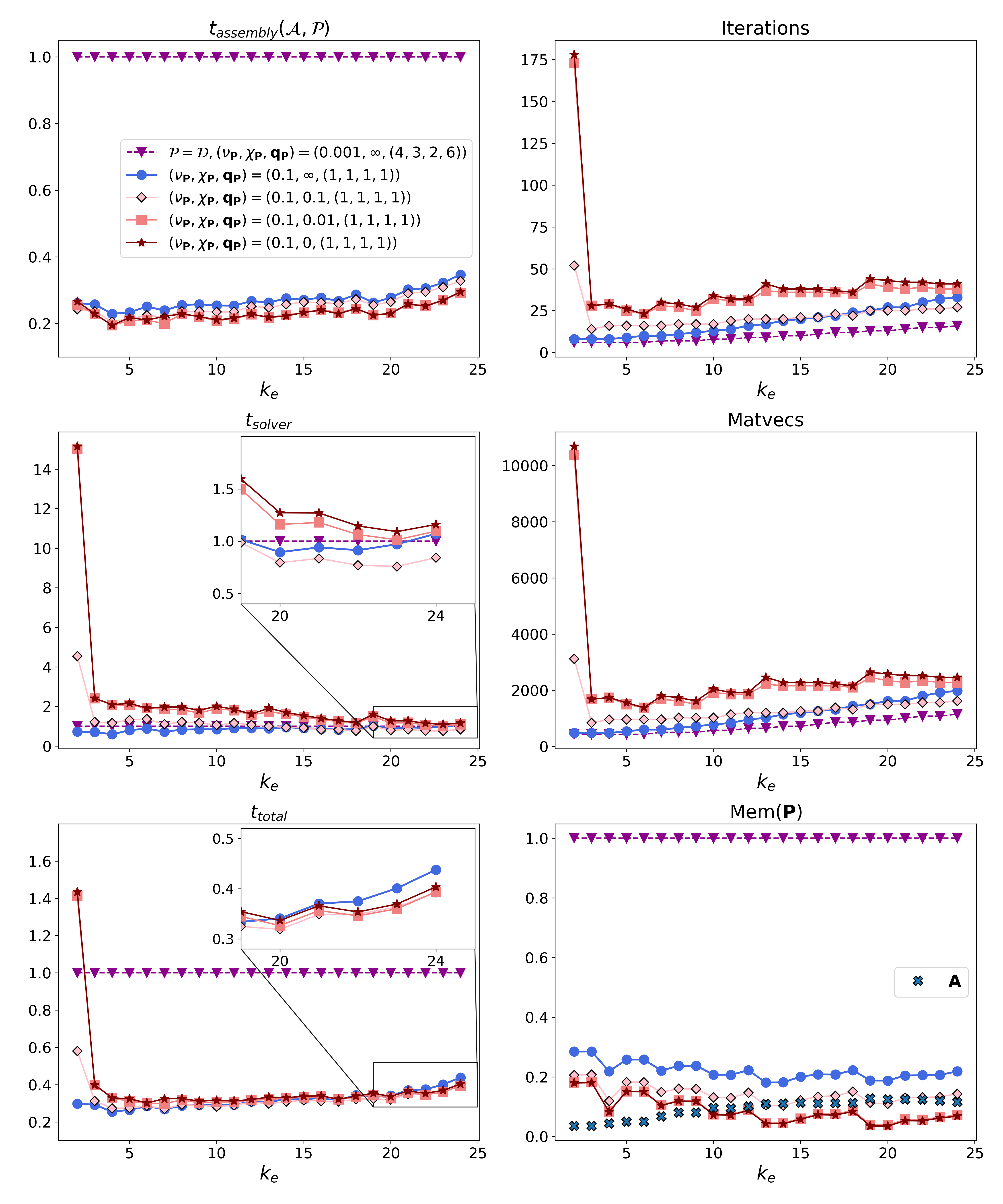}
    \caption{Analogue of Figure \ref{fig:biparametric_D_vs_k} for $\bfcal{P}=\bfcal{D}^i$. The reference solution is $\bfcal{P}=\bfcal{D}$ with $\nu_{\mathbf{P}}=\nu_{\mathbf{A}}=0.001$, $\chinf_{\mathbf{P}}=\chinf_{\mathbf{A}}=\infty$, $\mathbf{q}_{\mathbf{P}}=\mathbf{q}_{\mathbf{A}}=(4,3,2,6)$.}
\label{fig:biparametric_Di_vs_k}
\end{figure}

\begin{figure}[thp!]
    \centering
    \includegraphics[width = \textwidth]{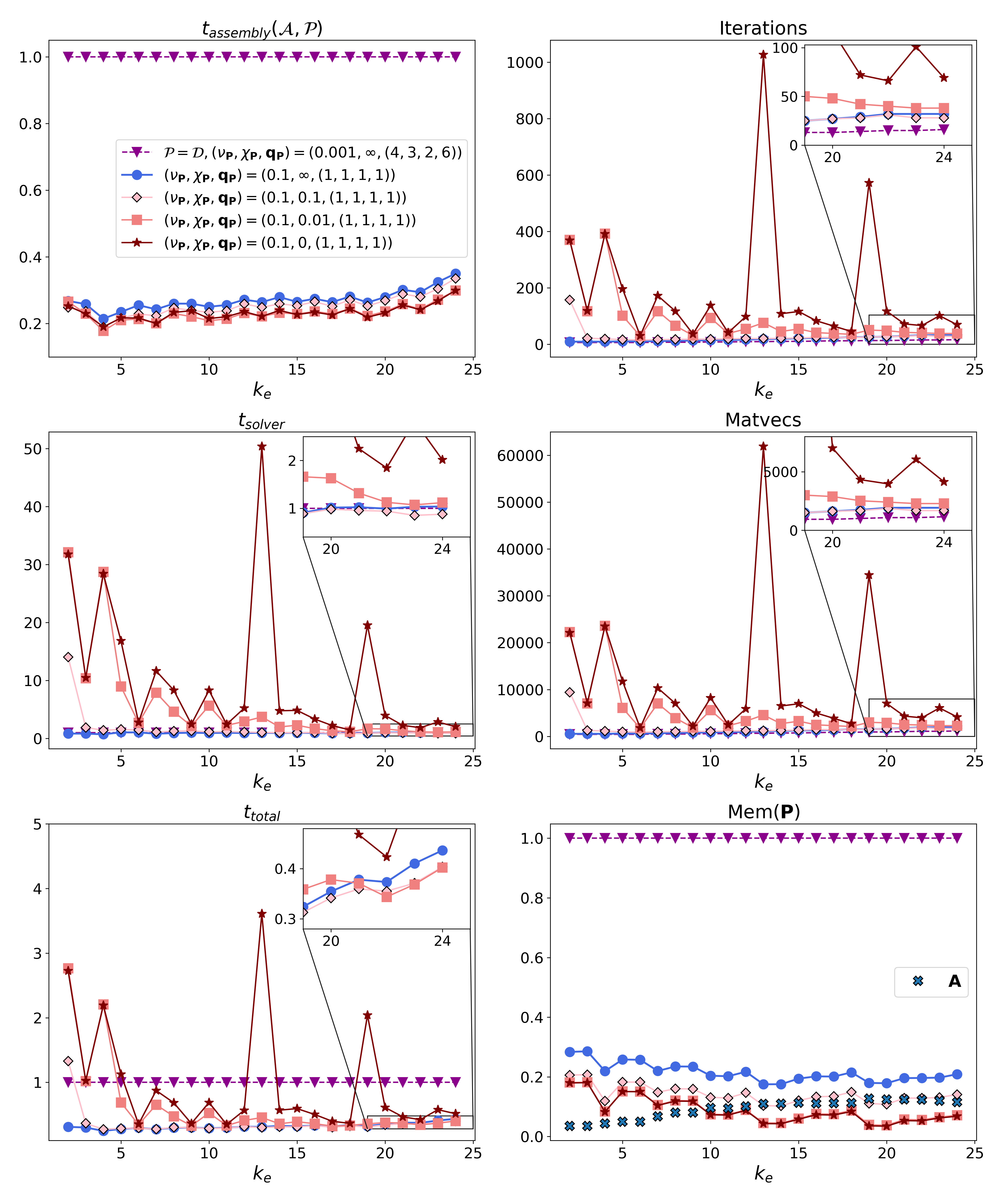}
    \caption{
    Analogue of Figure \ref{fig:biparametric_Di_vs_k} for $\bfcal{P}=\bfcal{D}^e$.} \label{fig:biparametric_De_vs_k}
\end{figure}

\begin{figure}[thp!]
    \centering
    \includegraphics[width = \textwidth]{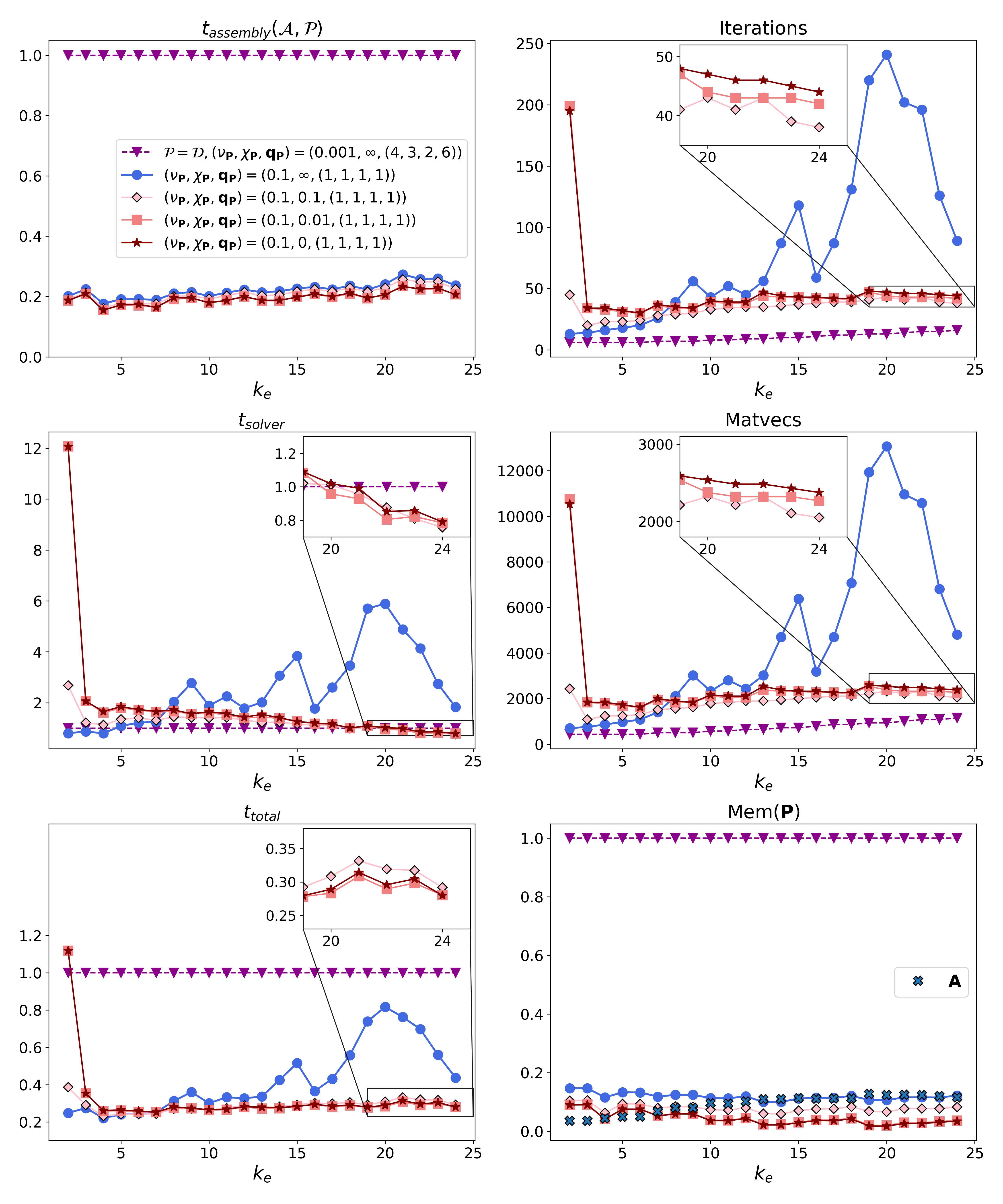}
    \caption{
Analogue of Figure \ref{fig:biparametric_Di_vs_k} for $\bfcal{P}=\bfcal{S}^i$.} \label{fig:biparametric_Si_vs_k}
\end{figure}

\begin{figure}[thp!]
    \centering
    \includegraphics[width = \textwidth]{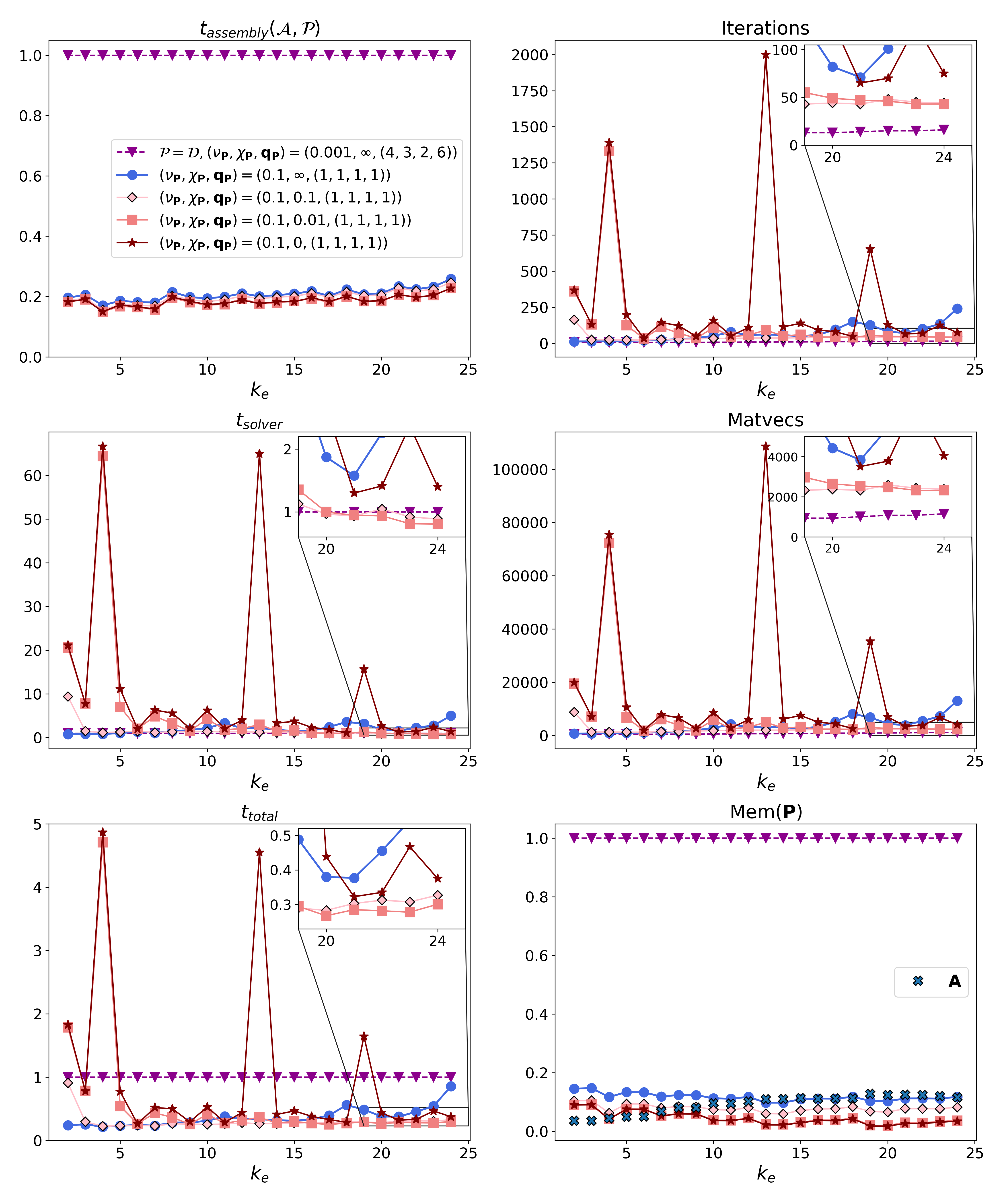}
    \caption{
Analogue of Figure \ref{fig:biparametric_Di_vs_k} for $\bfcal{P}=\bfcal{S}^e$.} \label{fig:biparametric_Se_vs_k}
\end{figure}

Analogous results for the reduced preconditioners $\bfcal{P} = \bfcal{D}^i, \ \bfcal{D}^e, \ \bfcal{S}^i$ and $\bfcal{S}^e$ are presented in in Figures \ref{fig:biparametric_Di_vs_k}--\ref{fig:biparametric_Se_vs_k}. 
As before, timings and memory costs are normalised relative to those for $\bfcal{P}=\bfcal{D}$ with $(\nu_{\mathbf{P}},\chinf_{\mathbf{P}},\mathbf{q}_{\mathbf{P}})=(\nu_{\mathbf{A}},\chinf_{\mathbf{A}},\mathbf{q}_{\mathbf{A}})=(0.001,\infty,(4,3,2,6))$. 
The results for $\bfcal{P} = \bfcal{D}^i$ in Figure \ref{fig:biparametric_Di_vs_k} are similar to those for $\bfcal{P} = \bfcal{D}$ in Figure \ref{fig:biparametric_D_vs_k}, in the sense that the significant increase in iteration/matvec count for small $\chi_{\mathbf{P}}$ is only observed for low-frequencies, but that increase is reduced compared to that for $\bfcal{P} = \bfcal{D}$. The total time is below 40\% (except for low-frequencies) while the memory cost drops below 30\% and becomes smaller compared to the memory cost of $\mathbf{A}$ at the higher frequencies, for $\chi_\mathbf{P} = 0.01$ and 0. The results for $\bfcal{P} = \bfcal{D}^e$ in Figure \ref{fig:biparametric_De_vs_k} are more erratic, however, with large spikes in iteration/matvec count (and hence solver time) for small values of $\chi_{\mathbf{P}}$ at certain values of $k_e$, indicating some kind of instability. At these values of $k_e$, while the memory cost is the same as that of $\bfcal{P} = \bfcal{D}^i$, total solution time for $\bfcal{P} = \bfcal{D}^e$ often exceeds that of the reference case.

The results for the reduced preconditioner $\bfcal{P} = \bfcal{S}^i$ in Figure \ref{fig:biparametric_Si_vs_k} show that for a bi-parametric implementation with $\chi_\mathbf{P} = \infty$ we still see the erratic behaviour at high $k_e$ observed for the non-bi-parametric implementation of Figure \ref{fig:reducedCalderon_multiple_vs_k}. However, reducing $\chi_\mathbf{P}$ to 0.1, 0.01 and 0 remedies this, with total time consistently between 20-40\% of the reference case (except at the lowest frequency), and the memory cost below 20\%, going as low as 2\% for higher $k_e$ and $\chi_\mathbf{P} = 0$. This makes 
a bi-parametric implementation of $\bfcal{P} = \bfcal{S}^i$ with $\chi_\mathbf{P}=0$ an attractive choice, one that we will use again in the next section. 
Finally, for $\bfcal{P} = \bfcal{S}^e$ (Figure \ref{fig:biparametric_Se_vs_k}), the memory cost is the same as that of $\bfcal{P} = \bfcal{S}^i$, but, like for $\bfcal{P} = \bfcal{D}^e$, we observe large iteration/matvec counts (and hence solver times) for small values of $\chi_{\mathbf{P}}$ at certain values of $k_e$, making the total time larger than that of the reference case. 

\section{Application: Electromagnetic scattering by ice crystals}
\label{sec:Applications}

Having validated our acceleration techniques on a simple benchmark problem, we now apply them to a large-scale problem relevant to atmospheric physics applications. Specifically, we consider electromagnetic scattering of an incident plane wave by the 8-branch ice crystal aggregate studied in \cite{yang1998single}, illustrated in Figure \ref{fig:geometry} (right panel). 
Such an aggregate is representative of ice crystal aggregates found in cirrus clouds, scattering from which is an important component in climate modelling \cite{baran2009review,baran2012single,liou2016light}.

\begin{table}[t!]
\centering
\begin{tabular}{ccccc}
\toprule
Frequency $f$ & Refractive index $n$ &  Wavelength $\lambda_e$ & Size parameter $\pi D_{max}/\lambda_e$ & \# dofs $N$   \\
\midrule
$\ 50$ GHz & $1.7754+0.00066i$ & 0.60 cm & 5  & 2556   \\[5pt]
$183$ GHz & $1.7754+0.00243i$ & 0.16 cm & 20 & 26418  \\[5pt]
$325$ GHz & $1.7754+0.00440i$ & 0.092 cm & 34 & 81318  \\[5pt]
$664$ GHz & $1.7754+0.00972i$ & 0.045 cm & 70 & 332523 \\[4pt]
\bottomrule
\end{tabular}
\caption{Refractive indices of ice at $-40^\circ$ C (from \cite{matzler2006microwave}), wavelengths, size parameters and number of dofs in the Bempp approximation space, at the four studied frequencies, for scattering by the 8-branch ice crystal aggregate in Figure \ref{fig:geometry} (right panel). }
\label{table:8aggr_info}
\end{table} 

We scale the aggregate so as to have 
diameter 
(maximum dimension) 
$D_{max} = 1$ cm, and modify the positions of the constituent crystals very slightly from those in \cite{yang1998single} so that they are non-overlapping, allowing the aggregate to be treated as union of $M=8$ disjoint scatterers. 
Although our previous experiments consisted of a simple set up of three cubes with significant distance between them, the 8-branch aggregate is a good example of compact ice crystals and provides a good test case of our accelerated methods in a multi-particle set up with minimal distance between the scatterers.
Slightly modifying the positions of the scatterers should not affect the accuracy of the solution or the resulting scattering properties of the ice crystals as the scales of the gaps are much smaller than the dimensions of the individual monomers or the resulting aggregate.
We note that small variations in the geometries of the scatterers in order to create valid BEM domains are not uncommon in the literature; \cite{groth2015boundary} for example introduced small cubes at the points where corners of individual monomers meet in order to ensure Lipschitz continuity of the resulting aggregate and treated it as a single particle problem.
In \cite{kleanthous2019calderon}, the same aggregate was treated as a multi-particle problem with the cube removed and the individual monomers touching at corners. The multi-particle preconditioned PMCHWT formulation still produced an accurate result and performed well.

We assume that the incident wave has the form $\left[0, \exp \left( \frac{ \mathrm{i} 2\pi f (x+z)}{\sqrt{2} c}\right), 0 \right] $, where $f$ is the frequency (in Hz) and $c$ is the speed of light. We focus our attention on four specific frequencies, $f=$ 50 GHz, 183 GHz, 325 GHz, and 664 GHz, which are typical frequencies used by the atmospheric physics community for the purposes of microwave remote sensing and for the assimilation of all-sky cloud radiances in numerical weather models \cite{fox2019airborne}. The corresponding approximate refractive indices of ice at $-40^\circ$ C (from \cite{matzler2006microwave}), exterior wavelengths $\lambda_e$, size parameters $\pi D_{max}/\lambda_e$ and the number $N$ of dofs in the Bempp approximation space can be found in Table \ref{table:8aggr_info}.







\begin{figure}[thp!]
    \centering
    \includegraphics[width = \textwidth]{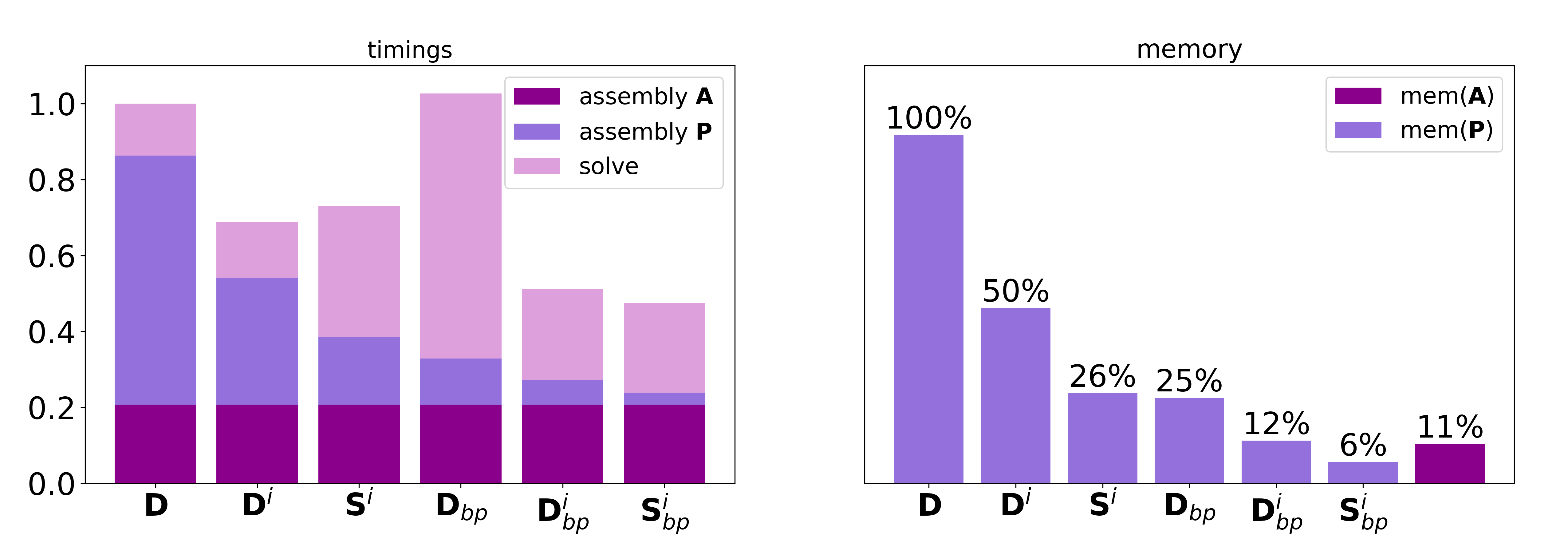} \\
    \includegraphics[width = 0.8 \textwidth]{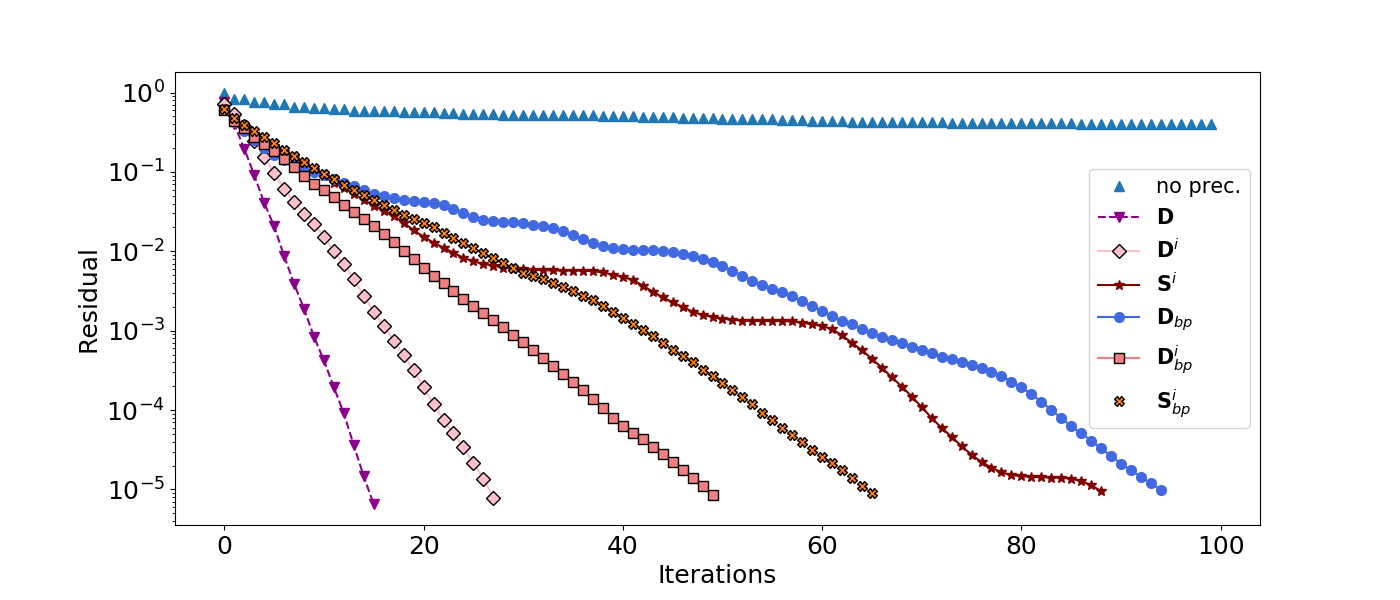}
    \caption{Performance at 50 GHz. Timings are normalised relative to the total time for $\bfcal{D}$, and memory costs are normalised relative to the memory cost of $\bfcal{D}$.
    As a reference, the GMRES iteration for the non-preconditioned operator, $\bfcal{A}$ (not shown in the top figures), was terminated after 1000 iterations (without having achieved the target tolerance), at which point it had already exceeded the total time taken for assembly and solution of $\bfcal{DA}$.}
    \label{fig:8aggregate_50GHz}
\end{figure}

\begin{figure}[thp!]
    \centering
    \includegraphics[width = \textwidth]{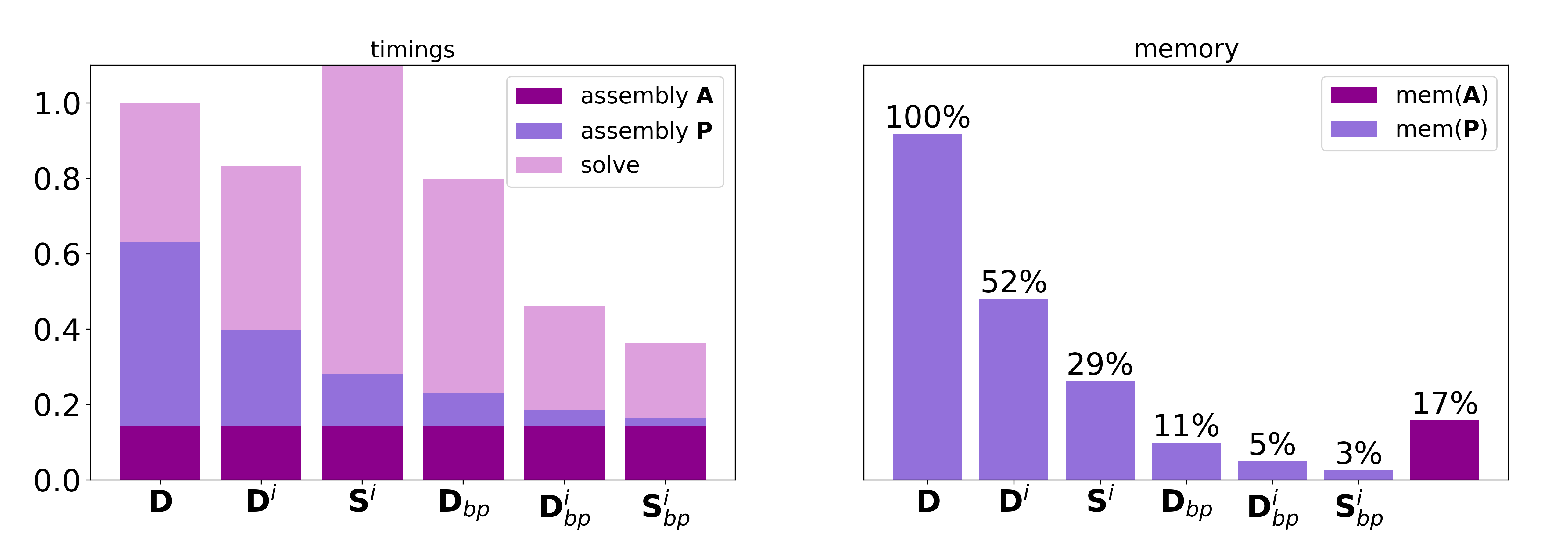} \\
    \includegraphics[width = 0.8 \textwidth]{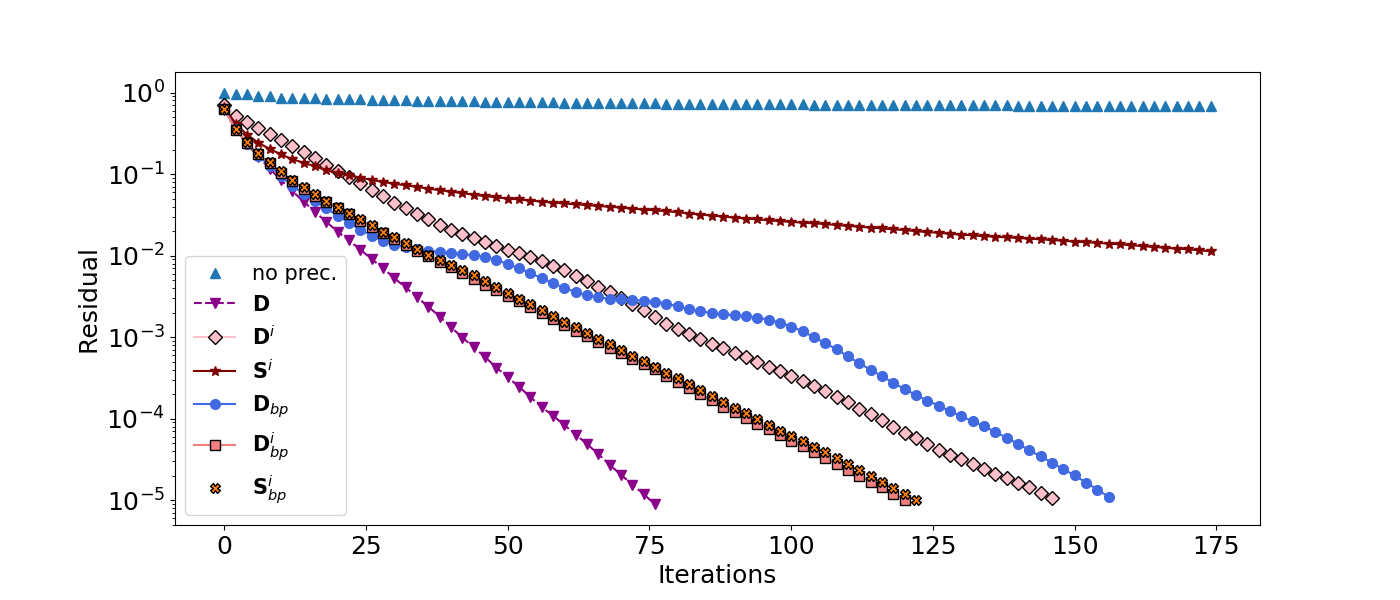}
    \caption{Performance at 183GHz. Timings are normalised relative to the total time for $\bfcal{D}$, and memory costs are normalised relative to the memory cost of $\bfcal{D}$.
    As a reference, the GMRES iteration for the non-preconditioned operator, $\bfcal{A}$ (not shown in the top figures), was terminated after 1500 iterations (without having achieved the target tolerance), at which point it had already exceeded the total time taken for assembly and solution of $\bfcal{DA}$.
    The GMRES iteration for the non-bi-parametric version of $\bfcal{S}^i$ (third from left in the graphs) was terminated after $400$ iterations (without having achieved the target tolerance). }
    \label{fig:8aggregate_183GHz}
\end{figure}

\begin{figure}[thp!]
    \centering
    \includegraphics[width = \textwidth]{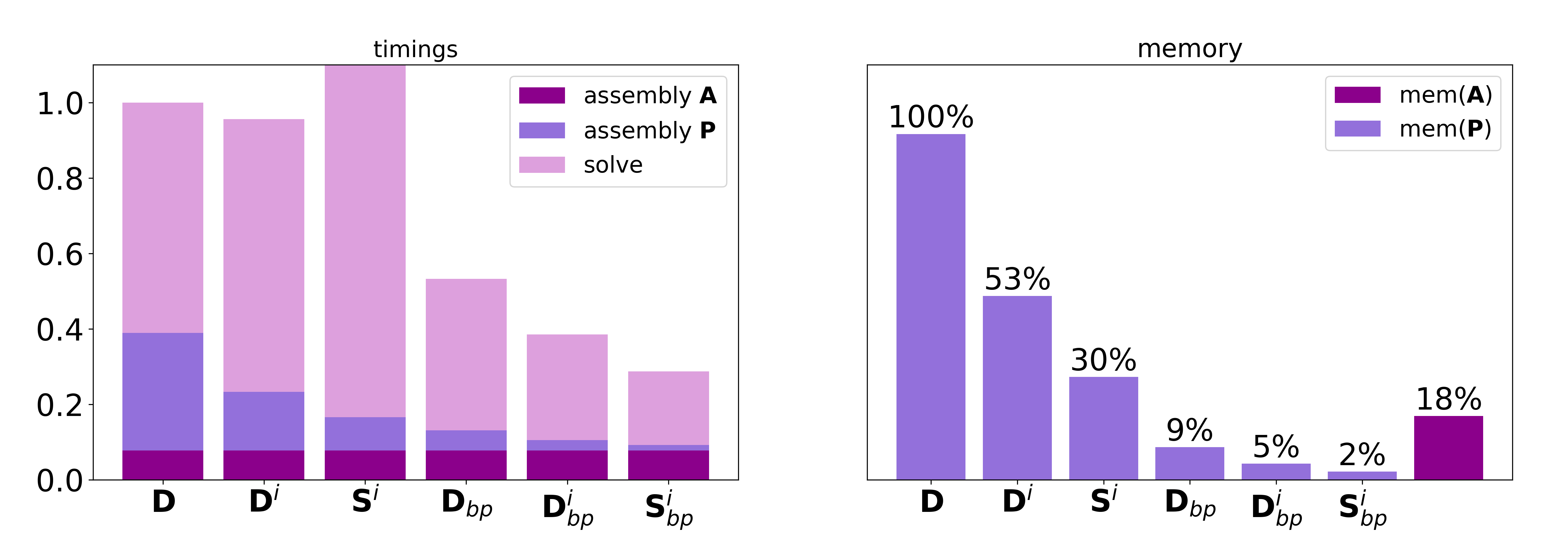}\\
    \includegraphics[width = 0.8 \textwidth]{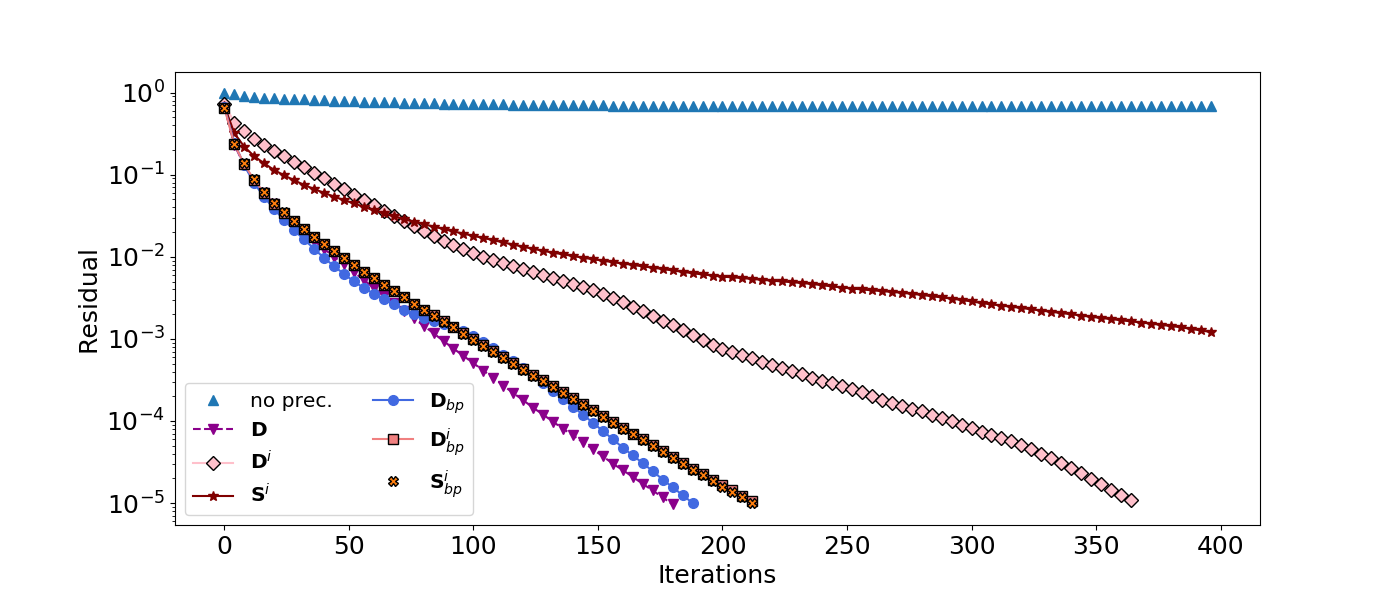}
    \caption{Performance at 325GHz. Timings are normalised relative to the total time for $\bfcal{D}$, and memory costs are normalised relative to the memory cost of $\bfcal{D}$. 
    The GMRES iterations for the non-preconditioned operator, $\bfcal{A}$ (not shown in the figure) and for the non-bi-parametric version of $\bfcal{S}^i$ (third from left in the graphs) were terminated after $800$ iterations (without having achieved the target tolerance). }
    \label{fig:8aggregate_325GHz}
\end{figure}

In Figures \ref{fig:8aggregate_50GHz}-\ref{fig:8aggregate_325GHz} we report computation times, memory costs and residual error curves for the frequencies 50 GHz, 183 GHz and 325 GHz respectively, for six choices of preconditioner: the block-diagonal preconditioner $\bfcal{D}$ and the reduced preconditioners $\bfcal{D}^i$ and $\bfcal{S}^i$, and their bi-parametric versions, labelled $\bfcal{D}_{bp}$, $\bfcal{D}^i_{bp}$ and $\bfcal{S}^i_{bp}$.
In all cases we assemble the operator $\bfcal{A}$ using the parameters $(\nu_{\mathbf{A}},\chinf_{\mathbf{A}},\mathbf{q}_{\mathbf{A}})=(0.001,\infty,(4,3,2,6))$. For the non-bi-parametric versions we take $(\nu_{\mathbf{P}},\chinf_{\mathbf{P}},\mathbf{q}_{\mathbf{P}})=(\nu_{\mathbf{A}},\chinf_{\mathbf{A}},\mathbf{q}_{\mathbf{A}})$ while for the bi-parametric versions we take $(\nu_{\mathbf{P}},\chinf_{\mathbf{P}},\mathbf{q}_{\mathbf{P}})=(0.1,0,(1,1,1,1))$. Our decision to set $\chi_{\mathbf{P}}=0$ was made to keep memory costs as low as possible, since our focus is on solving high frequency problems.  We do not present results for the exterior versions $\bfcal{D}^e$ and $\bfcal{S}^e$, since their behaviour in \S\ref{sct:benchmarks} was found to be erratic, leading to longer solver times compared to those of $\bfcal{D}^i$ and $\bfcal{S}^i$, respectively. 
GMRES parameters were kept the same as for previous experiments, i.e. $tol = 10^{-5}$ and $restart = 200$.

Our main observations from the results in Figures \ref{fig:8aggregate_50GHz}-\ref{fig:8aggregate_325GHz} are that, in line with the results in \S\ref{sct:benchmarks}, use of the reduced preconditioners significantly reduces assembly time and memory cost compared to the block-diagonal preconditioner $\bfcal{D}$. Furthermore, adopting a bi-parametric approach, with increased $\mathcal{H}$-matrix tolerance, reduced quadrature orders and far-field interactions neglected, brings the assembly time and memory cost well below that of $\bfcal{A}$. At the lowest frequency (50 GHz), solver time for the bi-parametric implementations is larger than that for $\bfcal{D}$, but, for the higher frequencies (183 GHz and 325 GHz) the reduced bi-parametric preconditioners $\bfcal{D}^i_{bp}$ and $\bfcal{S}^i_{bp}$ actually gave lower solver times than that for $\bfcal{D}$, the increased number of GMRES iterations (not reported here) being more than compensated for by the reduction in cost per iteration. 

\begin{figure}[thp!]
    \centering
    \includegraphics[width = 0.8 \textwidth]{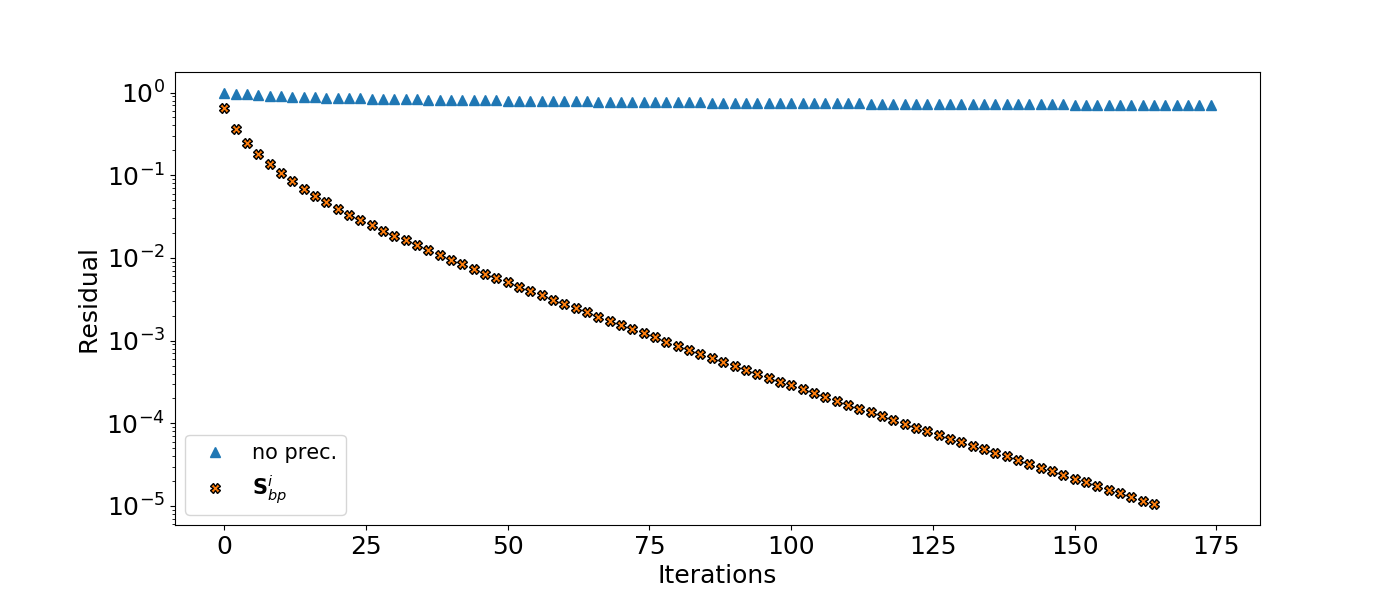}
    \caption{GMRES residual error at 664GHz.}
    \label{fig:8aggregate_664GHz}
\end{figure}

For 664 GHz, memory constraints meant that we were unable to assemble the matrices for any of the methods except for $\bfcal{S}^i_{bp}$. For this case, the assembly time for $\bfcal{S}^i_{bp}$ was 10 minutes, and 32 minutes for $\bfcal{A}$. GMRES converged in 62 minutes with 166 iterations.
GMRES residual curves for $\bfcal{S}^i_{bp}$ and the non-preconditioned case are shown in Figure \ref{fig:8aggregate_664GHz}.
The memory cost for $\bfcal{S}^i_{bp}$ was approximately 9 GB, just 8\% of the memory cost for the operator $\bfcal{A}$, which was 109 GB. Furthermore, extrapolating from the behaviour at 183 GHz and 325 GHz, the memory required for $\bfcal{S}^i_{bp}$ is about 1-2\% of the memory that would be required for the non-bi-parametric block-diagonal preconditioner $\bfcal{D}$, which is estimated to be approximately 650 GB with estimated assembly time at 160 minutes. Extrapolating from the solve time at 183 GHz and 325 GHz, we expect the total saving in computation time to be at least 75\%. 

For completeness, we repeated the experiments of this section using only the real part of the complex refractive indices in Table \ref{table:8aggr_info}. Although not reported here, our findings were similar to those reported in Figures \ref{fig:8aggregate_50GHz} -\ref{fig:8aggregate_325GHz} in terms of percentage time/memory savings, which suggests that the reduced preconditioners and bi-parametric implementation are also effective for cases of purely real interior and exterior wavenumbers.

\begin{figure}[t!]
    \centering
    \includegraphics[width = 0.48 \textwidth]{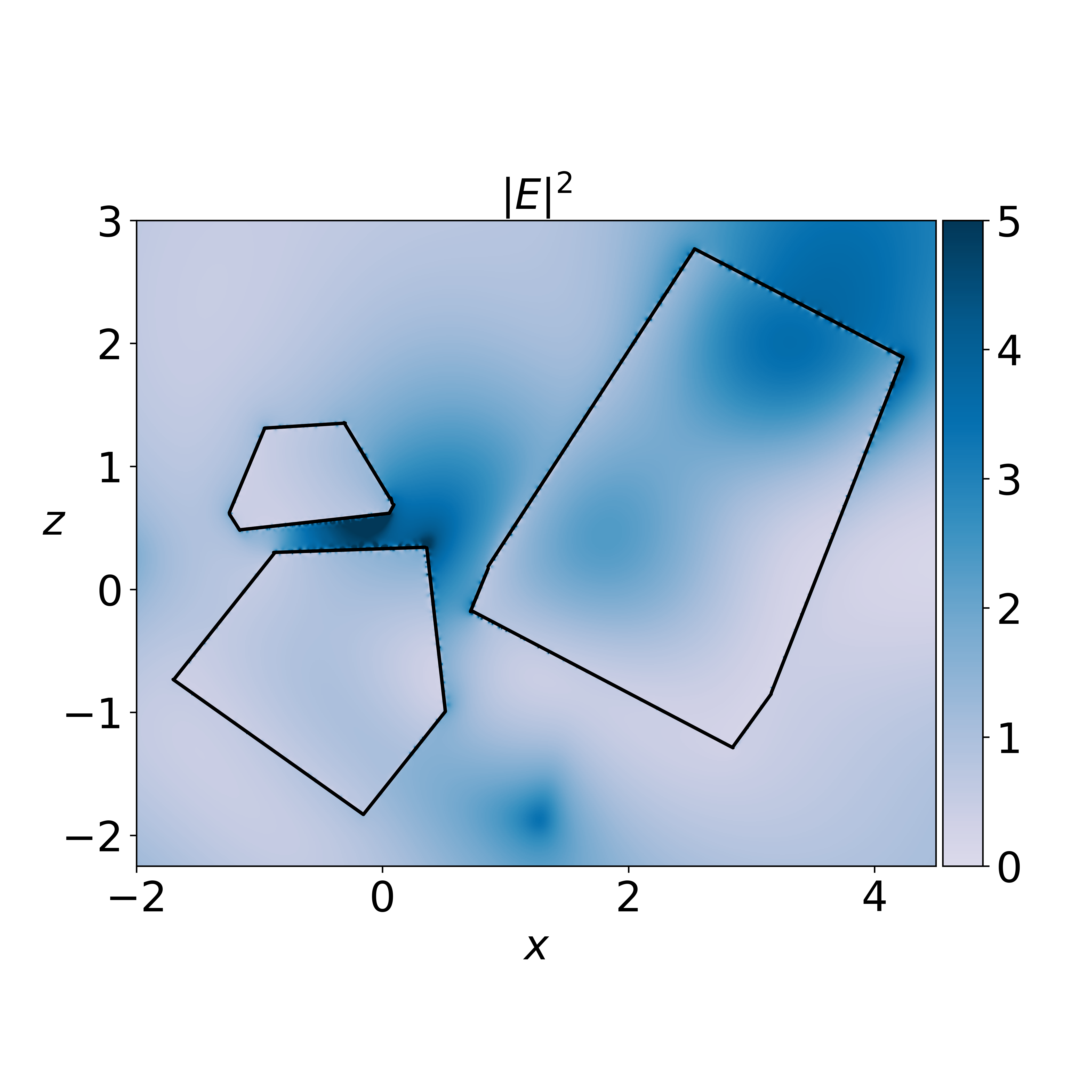}
    \hfill
    \includegraphics[width = 0.48 \textwidth]{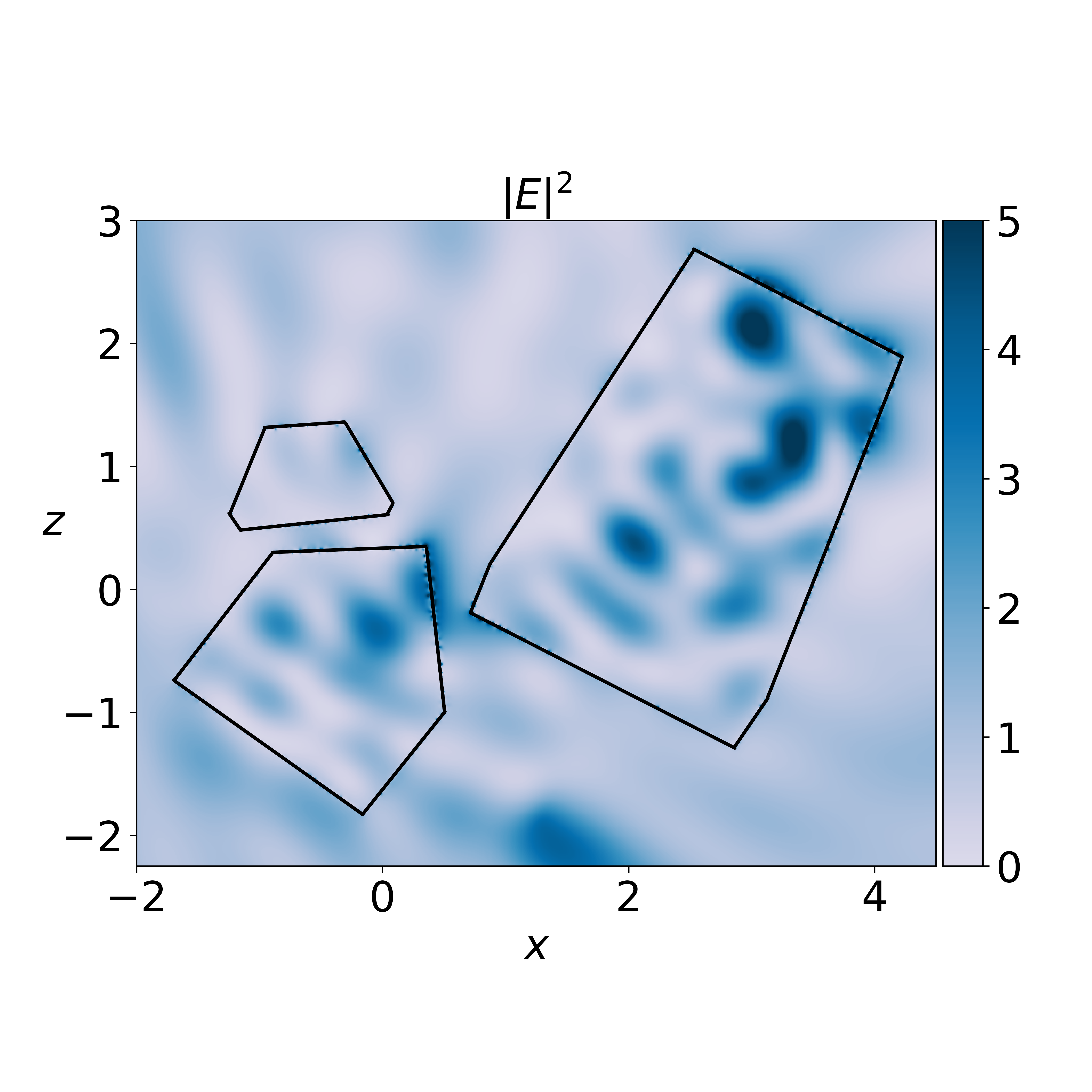}\\
    \includegraphics[width = 0.48 \textwidth]{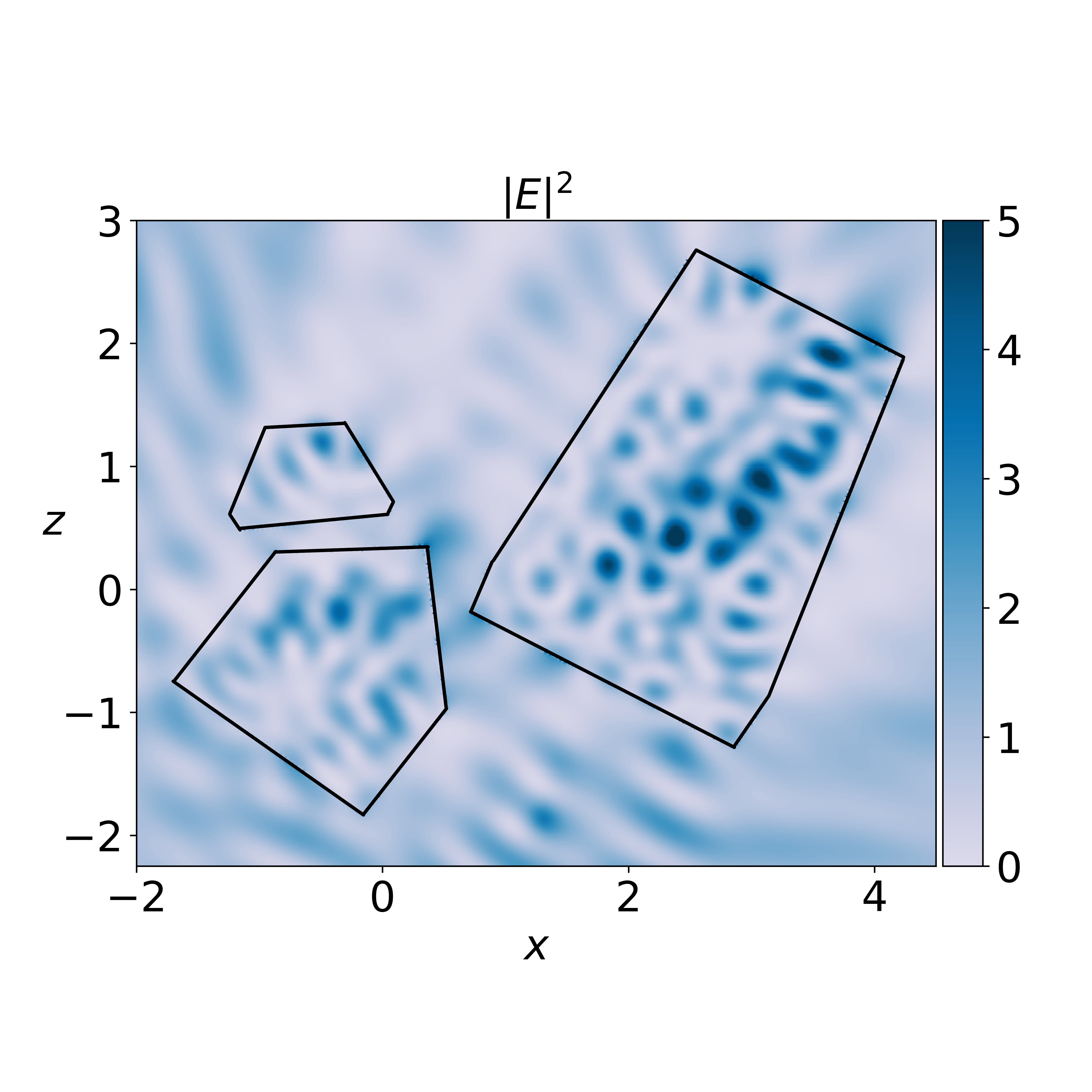}
    \hfill
	\includegraphics[width = 0.48 \textwidth]{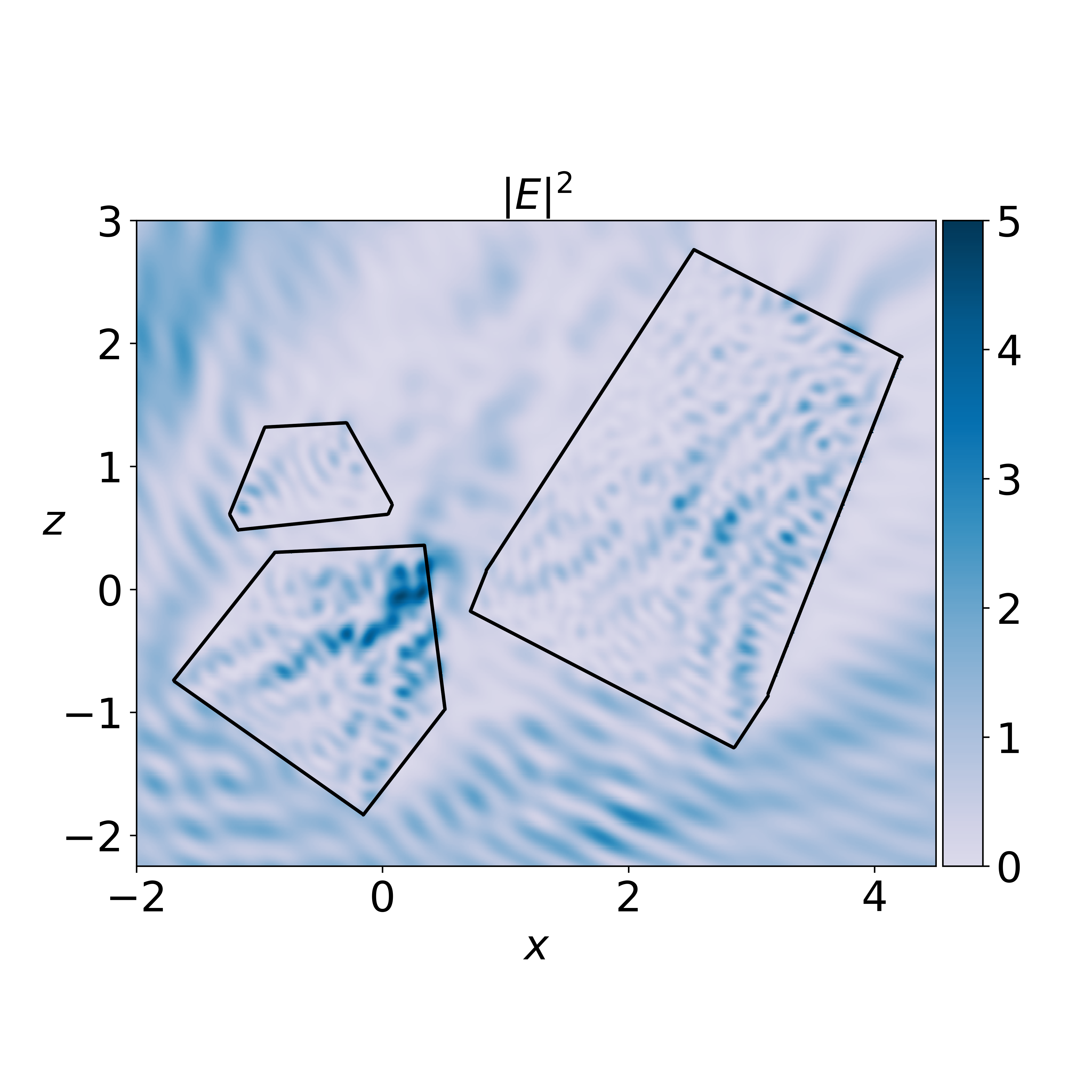}
    \caption{Square magnitude $|\mathbf{E}|^2$ of the electric field for scattering by the 8-branch aggregate of Figure \ref{fig:geometry} (right panel) in the plane $y=1$, at frequencies 50 GHz (top left), 183 GHz (top right), 325 GHz (bottom left) and 664 GHz (bottom right). Computations were done using the reduced bi-parametric preconditioner $\bfcal{P}=\bfcal{S}^i_{bp}$. }
    \label{fig:8aggr_plots}
\end{figure}

Finally, plots of the square magnitude $|\mathbf{E}|^2$ of the electric field in the plane $y=1$, for the four frequencies studied (and refractive indices as in Table \ref{table:8aggr_info}), computed using $\bfcal{P}=\bfcal{S}^i_{bp}$, are shown in Figure \ref{fig:8aggr_plots}. These plots clearly show the interior fields inside each of the constituent ice crystals increasing in complexity as the frequency increases, due to the interference between the multiply-scattered fields.

\section{Conclusion}\label{sct:conclusion}
In this paper, we investigated techniques for accelerating Calder\'on preconditioning for the PMCHWT boundary integral equation formulation of Maxwell transmission problems involving multiple scatterers. Our main result is that the high cost of the barycentric refinement necessary for stable discretisation of operator products can be mitigated completely by using a reduced preconditioner (discarding a subset of the operators appearing in the block-diagonal preconditioner of \cite{kleanthous2019calderon}, which we took as our reference method) combined with a bi-parametric implementation (using a lower-quality $\mathcal{H}$-matrix assembly for the preconditioner than for the original operator and discarding far-field interactions).
For the application of a large-scale problem representing an ice crystal aggregate (with individual monomers closer together), assembly times and memory cost for our reduced bi-parametric preconditioners (on the dual mesh) were significantly lower compared to those of the original operator (on the primal mesh).
The results presented in Sections \ref{sct:benchmarks} and \ref{sec:Applications} demonstrate that the proposed methods are efficient in both simple and complex setups and for a range of frequencies.

The optimal choice of preconditioner and assembly parameters depends on a number of factors including the geometry of the scatterer, the size of the scatterer relative to the incident wavelength, and the computational resources available.  Given the erratic behaviour of $\bfcal{D}^e$ and $\bfcal{S}^e$ we do not recommend their use until further research into their performance is undertaken. For the application in \S\ref{sec:Applications}, the best performance at all frequencies considered, in terms of both (total) computational time and memory cost, was given by a bi-parametric implementation of the reduced preconditioner $\bfcal{S}^i$, using ACA parameter $\nu_{\mathbf{P}}=0.1$, Bempp quadrature orders $(1,1,1,1)$ and all far-field interactions neglected ($\chi_{\mathbf{P}}=0$). However, in our study of the benchmark problem in \S\ref{sct:benchmarks} and for the case of 50GHz in \S\ref{sec:Applications}, this choice was found to give long solver times, especially at low frequencies. 
For this reason, if multiple GMRES solves are to be performed at the same time, for example in the case of orientational averaging, at low frequencies where memory constraints are not a concern but optimising GMRES time is important, then a sensible choice would be to use $\bfcal{P} = \bfcal{D}^i$ with a bi-parametric implementation using all far field interactions, i.e. $\nu_\mathbf{P} = 0.1$, $\chi_\mathbf{P}= \infty$ and $\mathbf{q}_\mathbf{P} = (1,1,1,1)$. 
For all other applications involving high frequencies where both memory and total computation time are expensive, then the optimal choice would be to use $\bfcal{P} = \bfcal{S}^i$, with $\nu_{\mathbf{P}}=0.1$, quadrature orders $(1,1,1,1)$ and $\chi_{\mathbf{P}}=0$.

We remark that although the focus of this paper and any results presented were for the PMCHWT formulation, the ideas of accelerating Calder\'on preconditioning presented in this paper can also be applied to other formulations such as the M\"uller formulation \cite{muller2013foundations}. 
Application of the theory in  \cite{ESCAPILINCHAUSPE2021220} to our proposed acceleration techniques could lead to controlled (i) approximation errors, and (ii) convergence bounds for GMRES.

We note that the PMCHWT BIE formulation for the multi-particle problem shown in Section \ref{sct:scattering_problem} can be extended to scattering by objects with inclusions or objects with multiple layers of homogeneous dielectric materials.
The former can be used to simulate scattering by ice crystals with air bubbles (such as in  \cite{yang2000efficient, labonnote2001polarized, xie2009effect, groth2015boundary}) or other trapped particles such as soot impurities \cite{yang2000efficient} or mineral aerosol \cite{labonnote2001polarized}.
The latter can be used for ice crystals falling through the atmosphere which are covered by a liquid surface due to the surface of the ice crystal melting (i.e. sticky ice).
The formulation for these cases has been briefly introduced in \cite{kleanthous2021accelerated} and the ideas of accelerated Calder\'on preconditioning can also be applied there.

\section*{Acknowledgements}
The work of the first author was supported by NERC and the UK Met Office (CASE PhD studentship to A.\ Kleanthous, grant NE/N008111/1). D.P. Hewett acknowledges support from EPSRC, grant EP/S01375X/1. P.~Escapil-Inchausp\'e and C.~Jerez-Hanckes thank the support of Fondecyt Regular 1171491.

\bibliography{mybibfile}
\end{document}